%% file: APZV2.tex
\documentclass[11pt,oneside,reqno]{amsart}

\usepackage[a4paper, total={450pt,675pt}]{geometry}
\usepackage[OT2,T1]{fontenc}
\usepackage[utf8]{inputenc}
\usepackage{amssymb,bm}

\usepackage[dvipsnames]{xcolor}
\usepackage[
    colorlinks=true,
    linkcolor=Maroon,
    citecolor=JungleGreen,
    urlcolor=NavyBlue]{hyperref}

\usepackage[capitalize,nameinlink]{cleveref}

\usepackage{enumitem}
\usepackage{dsfont}
\usepackage{subcaption}
\usepackage{centernot}

\usepackage{tikz}
\usetikzlibrary{arrows}
\usetikzlibrary{decorations.pathreplacing,angles,quotes}
\usepackage{pgfplots}
\pgfplotsset{compat=1.15}
\usepackage{mathrsfs}
\usetikzlibrary{arrows}

\makeatletter
\def\@tocline#1#2#3#4#5#6#7{\relax
  \ifnum #1>\c@tocdepth 
  \else
    \par \addpenalty\@secpenalty\addvspace{#2}%
    \begingroup \hyphenpenalty\@M
    \@ifempty{#4}{%
      \@tempdima\csname r@tocindent\number#1\endcsname\relax
    }{%
      \@tempdima#4\relax
    }%
    \parindent\z@ \leftskip#3\relax
    \advance\leftskip\@tempdima\relax
    \rightskip\@pnumwidth plus4em \parfillskip-\@pnumwidth
    #5\leavevmode\hskip-\@tempdima
      \ifcase #1
       \or\or \hskip 2em \or \hskip 2em \else \hskip 3em \fi%
      #6\nobreak\relax
    \dotfill\hbox to\@pnumwidth{\@tocpagenum{#7}}\par
    \nobreak
    \endgroup
  \fi} 
\makeatother

\usepackage[
    backend=bibtex,
    style=alphabetic,
    maxbibnames=10,
    sorting=nyt]{biblatex}

\DeclareLabelalphaTemplate{
  \labelelement{
    \field[final]{shorthand}
    \field{label}
    \field[strwidth=3]{labelname}}
  \labelelement{
    \field[strwidth=2,strside=right]{year}}
}

\DeclareFieldFormat[article,book,incollection,misc]{title}{#1}
\DeclareFieldFormat[inbook,incollection]{booktitle}{#1}
\DeclareFieldFormat[misc]{date}{\textit{preprint} (#1)}
\DeclareFieldFormat{pages}{#1}
\DeclareFieldFormat{doi}{
\href{http://www.ams.org/mathscinet-getitem?mr=#1}{#1}}

\renewbibmacro{in:}{%
  \ifentrytype{article}{}{\printtext{\bibstring{in}\intitlepunct}}}

\newtheorem{theorem}{Theorem}[section]

\newtheorem{lemma}[theorem]{Lemma}
\newtheorem{proposition}[theorem]{Proposition}
\newtheorem{corollary}[theorem]{Corollary}

\newtheorem{remark}[theorem]{Remark}

\crefname{section}{Sect.}{section}
\numberwithin{equation}{section}

\newcommand*\diff{\mathop{}\!\mathrm{d}}
\DeclareMathOperator{\Ad}{Ad}

\DeclareMathOperator{\im}{Im}

\DeclareMathOperator{\supp}{supp}

\DeclareMathOperator{\const}{const.}
\DeclareMathOperator{\co}{co}

\begin{document}
\title[Asymptotic behavior for the heat equation on noncompact 
symmetric spaces]{Asymptotic behavior of solutions to the heat equation\\
on noncompact symmetric spaces}

\author{Jean-Philippe ANKER, Effie Papageorgiou and Hong-Wei ZHANG}

\begin{abstract}
This paper is twofold. The first part aims to study the long-time asymptotic
behavior of solutions
to the heat equation on Riemannian symmetric spaces 
$G/K$ of noncompact type and of general rank.
We show that any
solution to the heat equation with bi-$K$-invariant $L^{1}$ 
initial data behaves asymptotically as the mass times the fundamental solution,
and provide a counterexample 
in the non bi-$K$-invariant case.
These answer problems recently raised by J.L. Vázquez.
In the second part, we investigate the long-time 
asymptotic behavior of solutions to the heat equation associated 
with the so-called distinguished Laplacian on $G/K$.
Interestingly, we observe in this
case phenomena which are similar to the
Euclidean setting, namely $L^1$ asymptotic convergence with no bi-$K$-invariance
condition and strong $L^{\infty}$ convergence.
\end{abstract}


\makeatletter
\@namedef{subjclassname@2020}{\textnormal{2020}
    \it{Mathematics Subject Classification}}
\makeatother

\maketitle
\tableofcontents
\section{Introduction}\label{Section.1 Intro}

The heat equation is one of the most fundamental partial differential equations 
in mathematics. 
After the celebrated work of {\it Joseph Fourier} in 1822,
it has been extensively studied in various settings and is known to play a
central role in several areas of mathematics (see for instance \cite{Gri2009}).
The following classical long-time asymptotic convergence result
corresponds to the {\it Central Limit Theorem} of probability in the PDE setting. 
We refer to the expository survey \cite{Vaz2018} for more details on this property.

\begin{theorem}
Consider the heat equation
\begin{align}\label{S1 HE intro}
\begin{cases}
    \partial_{t}u(t,x)\,
    =\,\Delta_{\mathbb{R}^{n}}u(t,x),
    \qquad\,t>0,\,\,x\in\mathbb{R}^{n}\\[5pt]
    u(0,x)\,=\,u_{0}(x),
\end{cases}
\end{align}
where the initial data $u_{0}$ belongs to $L^{1}(\mathbb{R}^n)$.
Denote by $M=\int_{\mathbb{R}^n}\diff{x}\,u_{0}(x)$ the mass of $u_{0}$
and by $G_{t}(x)\,=\,(4\pi{t})^{-n/2}e^{-|x|^{2}/4t}$ the heat kernel.
Then the solution to \eqref{S1 HE intro} satisfies:
\begin{align}
    \|u(t,\,\cdot\,)-MG_{t}\|_{L^{1}(\mathbb{R}^n)}\,
    \longrightarrow\,0
\end{align}
and
\begin{align}\label{S1 Linf R}
    t^{\frac{n}{2}}\,\|u(t,\,\cdot\,)-MG_{t}\|_{L^{\infty}(\mathbb{R}^n)}\,
    \longrightarrow\,0
\end{align}
as $t\rightarrow\infty$. The $L^p$ ($1<p<\infty$) norm estimates 
follow by convexity.
\end{theorem}

To our knowledge, analogous properties on negatively curved manifolds
were first investigated by Vázquez in his recent work \cite{Vaz2019},
which deals with real hyperbolic spaces and where he suggests possible generalizations,
in particular to all Riemannian symmetric spaces of noncompact type
(other hyperbolic spaces in rank one and noncompact symmetric spaces of higher rank).
This is fully achieved in the present paper,
where we also clarify some arguments in \cite{Vaz2019}
and consider a related Laplacian with Euclidean type properties. Let us elaborate.
Given a symmetric space $\mathbb{X}=G/K$ of non compact type,
let us denote by $\Delta$ the Laplace-Beltrami operator on $\mathbb{X}$
and by $h_t$ the associated heat kernel, i.e.,
the bi-$K$-invariant convolution kernel of the semi-group $e^{t\Delta}$.
Our first main result is about the long-time
convergence of solutions to the heat equation:
\begin{align}\label{S1 HE X}
    \partial_{t}u(t,x)\,=\,\Delta_{x}u(t,x),
    \qquad
    u(0,x)\,=\,u_{0}(x).
\end{align}

\begin{theorem}\label{S1 Main thm 1}
Let $u_{0}\in{L^{1}(\mathbb{X})}$ be a bi-$K$-invariant initial data
and $M=\int_{\mathbb{X}}\diff{x}\,u_{0}(x)$ be its mass. Then the solution to
the heat equation \eqref{S1 HE X} satisfies
\begin{align}\label{S1 Main thm 1 convergence}
    \|u(t,\,\cdot\,)-Mh_{t}\|_{L^{1}(\mathbb{X})}\,
    \longrightarrow\,0
    \qquad\textnormal{as}\quad\,t\rightarrow\infty.
\end{align}
Moreover, this convergence fails in general without the bi-$K$-invariance assumption.
\end{theorem}

\begin{remark}
The convergence \eqref{S1 Main thm 1 convergence} generalizes the result 
previously obtained in \cite{Vaz2019} on $\mathbb{H}^{n}(\mathbb{R})$.
Vázquez also provided a counterexample in the non bi-$K$-invariant case
on $\mathbb{H}^{3}(\mathbb{R})$, where the heat kernel $h_{t}$ has an 
explicit expression, and asked if it could be extended in other dimensions.
We establish such a counterexample not only for $\mathbb{H}^{n}(\mathbb{R})$,
$n\neq3$, but also for the symmetric spaces of general rank, see \cref{S3 Sub4}. Finally, our method sheds light to why this non-euclidean discrepancy occurs, see \cref{Remark 3.10}.
\end{remark}

\begin{remark}
If the initial data is in addition compactly supported, 
we obtain the better estimate
\begin{align*}
     \|u(t,\,\cdot\,)-Mh_{t}\|_{L^{1}(\mathbb{X})}\,
     \le\,C\,t^{-\frac12+\varepsilon}
     \qquad\forall\,t\ge1,
\end{align*}
where $C>0$ is a constant and $\varepsilon$ is any
small positive constant, see \cref{S3 Sub1} and 
\cref{S3 Sub2}.
\end{remark} 

\begin{remark}
We also provide the following sup norm (for which no bi-$K$-invariance is needed)
and $L^{p}$ $(1<p<\infty$) norm estimates:
\begin{align}
    \|u(t,\,\cdot\,)-Mh_{t}\|_{L^{\infty}(\mathbb{X})}\,
    &=\,
    \mathrm{O}\big(t^{-\frac{\nu}{2}}e^{-|\rho|^{2}t}\big)\label{LinftyX}\\[5pt]
    \|u(t,\,\cdot\,)-Mh_{t}\|_{L^{p}(\mathbb{X})}\,
    &=\,\mathrm{o}\big(t^{-\frac{\nu}{2p'}}
        e^{-\tfrac{|\rho|^{2}t}{p'}}\big)
\end{align}
as $t\rightarrow\infty$. 
Here, $p'$ denotes the dual exponent of $p$, defined by the formula
$\tfrac{1}{p}+\tfrac{1}{p'}=1$.
As previously observed on $\mathbb{H}^{n}(\mathbb{R})$,
the sup norm estimate in the present context is relatively weaker compared to
\eqref{S1 Linf R} in the Euclidean setting, while the $L^{p}$ norm estimate is 
similar.
Here, $\nu$ denotes the so-called dimension at infinity 
of $\mathbb{X}$ and $\rho$ is the 
half sum of positive roots with multiplicities, see \cref{Section.2 Prelim}.
\end{remark}

Let $S=N(\exp{\mathfrak{a}})=(\exp{\mathfrak{a}})N$ be the solvable group
occurring in the Iwasawa decomposition $G=N(\exp{\mathfrak{a}})K$. 
Then $S$ is identifiable, as a manifold, with the symmetric space 
$\mathbb{X}=G/K$.
Our second main contribution is to study the asymptotic convergence 
for solutions to the heat equation associated with the 
so-called distinguished Laplacian
$\widetilde{\Delta}$ on $S$. 
In order to state the results, let us introduce some indispensable notation,
which will be clarified in \cref{Section.2 Prelim}
and \cref{Section.4 Distinguished}.
Let $\mathfrak{a}$ be the Cartan subspace of $\mathbb{X}$. Denote by
$\varphi_{\lambda}$, where $\lambda\in\mathfrak{a}$, 
the spherical function occurring in the Harish-Chandra transform,
by $\widetilde{\delta}$ the modular function on $S$, 
and by $\widetilde{h}_{t}=\widetilde{\delta}^{1/2}e^{|\rho|^{2}t}h_{t}$
the fundamental solution to the Cauchy problem
\begin{align}\label{S1 HE S}
    \partial_{t}\widetilde{v}(t,g)\,
    =\,\widetilde{\Delta}_{g}\widetilde{v}(t,g),
    \qquad
    \widetilde{v}(0,g)\,=\,\widetilde{v}_{0}(g).
\end{align}
Let $\widetilde{\varphi}_{\lambda}=\widetilde{\delta}^{1/2}\varphi_{\lambda}$ 
be the modified spherical function and denote by 
$\widetilde{M}=\tfrac{\widetilde{v}_{0}*
\widetilde{\varphi}_{\lambda}}{\widetilde{\varphi}_{\lambda}}$ the mass function
on $S$ which generalizes the mass in the Euclidean case (see \cref{S4 mass remark}).
Then, we show the following long-time asymptotic convergence results.
\begin{theorem}\label{S1 Main thm 2}
Let $\widetilde{v}_{0}$ be a continuous and compactly supported initial data on
$S$. Then, the solution to the heat equation \eqref{S1 HE S} satisfies
\begin{align}\label{S1 L1 disting}
    \|\widetilde{v}(t,\,\cdot\,)-
    \widetilde{M}\,\widetilde{h}_{t}\|_{L^{1}(S)}\,
    \longrightarrow\,0
\end{align}
and
\begin{align}\label{S1 Linf disting}
    t^{\frac{\ell+|\Sigma_{r}^{+}|}{2}}
    \|\widetilde{v}(t,\,\cdot\,)-
    \widetilde{M}\,\widetilde{h}_{t}\|_{L^{\infty}(S)}\,
    \longrightarrow\,0
\end{align}
as $t\rightarrow\infty$. Here $\ell$ denotes the rank of $G/K$ 
and $\Sigma_{r}^{+}$ the set of positive reduced roots.
Analogous $L^p$ ($1<p<\infty$) norm estimates 
follow by convexity.
\end{theorem}

\begin{remark}
Let us comment about \eqref{S1 L1 disting} and \eqref{S1 Linf disting}.
Firstly, notice that the $L^1$ convergence \eqref{S1 L1 disting} holds 
with no bi-$K$-invariance restriction, in contrast to \cref{S1 Main thm 1}, 
and the sup norm estimate $\eqref{S1 Linf disting}$ is stronger than
(\ref{LinftyX}), as in the Euclidean
setting.
Secondly, the mass $\widetilde{M}$ is a bounded function and
no more a constant in general.
Thirdly, the power $\ell+|\Sigma_{r}^{+}|$,
which occurs in the time factor, is always different from the dimension at infinity
$\nu=\ell+2|\Sigma_{r}^{+}|$ and it is equal to the manifold dimension
$n=\ell+\sum_{\alpha\in\Sigma^{+}}m_{\alpha}$ 
if and only if the following equivalent conditions hold:
\begin{itemize}
    \item 
        the root system $\Sigma$ is reduced and all roots have multiplicity
        $m_{\alpha}=1$.
   \item
        $G$ is a normal real form.
\end{itemize}
\end{remark}

\begin{remark}
The asymptotic convergences \eqref{S1 L1 disting} and \eqref{S1 Linf disting}
hold for some larger classes of initial data, see \cref{Subsect other data}.
These classes are not optimal and finding the right function space is
an interesting question for further study.
\end{remark}

This paper is organized as follows. After the present introduction
in \cref{Section.1 Intro} and preliminaries in \cref{Section.2 Prelim}, 
we deal with the long-time asymptotic behavior of solutions to the heat equation
associated with the Laplace-Beltrami operator on symmetric spaces in 
\cref{Section.3 Asymp}. We start with smooth bi-$K$-invariant compactly supported
initial data, and prove, on the one hand, the long-time convergence in $L^{1}$ 
in the critical region where the heat kernel concentrates.
On the other hand, we show that both the solution and the heat kernel vanish
asymptotically outside that critical region. In the rest of this section,
we discuss these problems for more general initial data in the 
$L^{p}\,(p\ge1)$ setting, and provide a counterexample 
for $L^1$ in the non bi-$K$-invariant case.
In \cref{Section.4 Distinguished}, we investigate the asymptotic behavior
of solutions to the heat equation associated with the distinguished Laplacian.
After specifying the critical region in this context, we study the long-time
convergence in $L^1$ and in $L^{\infty}$ with compactly supported initial
data. Questions associated with other initial data are discussed at the end 
of the paper.

Throughout this paper, the notation
$A\lesssim{B}$ between two positive expressions means that 
there is a constant $C>0$ such that $A\le{C}B$. 
The notation $A\asymp{B}$ means that $A\lesssim{B}$ and $B\lesssim{A}$.

\section{Preliminaries}\label{Section.2 Prelim}
In this section, we first review spherical Fourier analysis 
on Riemannian symmetric spaces of noncompact type. The notation is standard 
and follows \cite{Hel1978,Hel2000,GaVa1988}.
Next we recall the asymptotic concentration of the heat kernel.
We refer to \cite{AnJi1999,AnOs2003} for more details on the heat kernel 
analysis in this setting.
\subsection{Noncompact Riemannian symmetric spaces}
Let $G$ be a semi-simple Lie group, connected, noncompact, with finite center, 
and $K$ be a maximal compact subgroup of $G$. The homogeneous space 
$\mathbb{X}=G/K$ is a Riemannian symmetric space of noncompact type.
Let $\mathfrak{g}=\mathfrak{k}\oplus\mathfrak{p}$ be the Cartan decomposition 
of the Lie algebra of $G$. The Killing form of $\mathfrak{g}$ induces 
a $K$-invariant inner product $\langle\,.\,,\,.\,\rangle$ on $\mathfrak{p}$, 
hence a $G$-invariant Riemannian metric on $G/K$.
We denote by $d(\,.\,,\,.\,)$ the Riemannian distance on $\mathbb{X}$.

Fix a maximal abelian subspace $\mathfrak{a}$ in $\mathfrak{p}$. 
The rank of $\mathbb{X}$ is the dimension $\ell$ of $\mathfrak{a}$.
We identify $\mathfrak{a}$ with its dual $\mathfrak{a}^{*}$ 
by means of the inner product inherited from $\mathfrak{p}$.
Let $\Sigma\subset\mathfrak{a}$ be the root system of 
$(\mathfrak{g},\mathfrak{a})$ and denote by $W$ the Weyl group 
associated with $\Sigma$. 
Once a positive Weyl chamber $\mathfrak{a}^{+}\subset\mathfrak{a}$ 
has been selected, $\Sigma^{+}$ (resp. $\Sigma_{r}^{+}$ 
or $\Sigma_{s}^{+}$)  denotes the corresponding set of positive roots 
(resp. positive reduced, i.e., indivisible roots or simple roots).
Let $n$ be the dimension and $\nu$ be the pseudo-dimension 
(or dimension at infinity) of $\mathbb{X}$: 
\begin{align}\label{S2 Dimensions}
\textstyle
    n\,=\,
    \ell+\sum_{\alpha \in \Sigma^{+}}\,m_{\alpha}
    \qquad\textnormal{and}\qquad
    \nu\,=\,\ell+2|\Sigma_{r}^{+}|
\end{align}
where $m_{\alpha}$ denotes the dimension of the positive root subspace
\begin{align*}
    \mathfrak{g}_{\alpha}\,
    =\,\lbrace{
    X\in\mathfrak{g}\,|\,[H,X]=\langle{\alpha,H}\rangle{X},\,
    \forall\,H\in\mathfrak{a}
    }\rbrace.
\end{align*}

Let $\mathfrak{n}$ be the nilpotent Lie subalgebra 
of $\mathfrak{g}$ associated with $\Sigma^{+}$ 
and let $N = \exp \mathfrak{n}$ be the corresponding 
Lie subgroup of $G$. We have the decompositions 
\begin{align*}
    \begin{cases}
        \,G\,=\,N\,(\exp\mathfrak{a})\,K 
        \qquad&\textnormal{(Iwasawa)}, \\[5pt]
        \,G\,=\,K\,(\exp\overline{\mathfrak{a}^{+}})\,K
        \qquad&\textnormal{(Cartan)}.
    \end{cases}
\end{align*}
Denote by $A(x)\in\mathfrak{a}$ and $x^{+}\in\overline{\mathfrak{a}^{+}}$
the middle components of $x\in{G}$ in these two decompositions, and by
$|x|=|x^{+}|$ the distance to the origin. For all $x,y\in{G}$, we have
\begin{align}\label{S2 Distance}
    |A(xK)|\,\le\,|x|
    \qquad\textnormal{and}\qquad
    |x^{+}-y^{+}|\,\le\,d(xK,yK),
\end{align}
see for instance \cite[Lemma 2.1.2]{AnJi1999}.
In the Cartan decomposition, the Haar measure 
on $G$ writes
\begin{align*}
    \int_{G}\diff{x}\,f(x)\,
    =\,
    |K/\mathbb{M}|\,\int_{K}\diff{k_1}\,
    \int_{\mathfrak{a}^{+}}\diff{x^{+}}\,\delta(x^{+})\, 
    \int_{K}\diff{k_2}\,f(k_{1}(\exp x^{+})k_{2})\,,
\end{align*}
with density
\begin{align}\label{S2 estimate of delta}
    \delta(x^{+})\,
    =\,\prod_{\alpha\in\Sigma^{+}}\,
        (\sinh\langle{\alpha,x^{+}}\rangle)^{m_{\alpha}}\,
    \asymp\,
        \prod_{\alpha\in\Sigma^{+}}
        \Big( 
        \frac{\langle\alpha,x^{+}\rangle}
        {1+\langle\alpha,x^{+}\rangle}
        \Big)^{m_{\alpha}}\,
        e^{2\langle\rho,x^{+}\rangle}
    \qquad\forall\,x^{+}\in\overline{\mathfrak{a}^{+}}. 
\end{align}
Here $K$ is equipped with its normalized Haar measure,
$\mathbb{M}$ denotes the centralizer of $\exp\mathfrak{a}$ in $K$ and the volume 
of $K/\mathbb{M}$ can be computed explicitly, see \cite[Eq (2.2.4)]{AnJi1999}.
Recall that $\rho\in\mathfrak{a}^{+}$ denotes the half sum of all positive roots 
$\alpha \in \Sigma^{+}$ counted with their multiplicities $m_{\alpha}$:
\begin{align*}
    \rho\,=\,
    \frac{1}{2}\,\sum_{\alpha\in\Sigma^{+}} \,m_{\alpha}\,\alpha.
\end{align*}

\subsection{Spherical Fourier analysis}
Let $\mathcal{S}(K \backslash{G}/K)$ be the Schwartz space of bi-$K$-invariant
functions on $G$. The spherical Fourier transform (Harish-Chandra transform)
$\mathcal{H}$ is defined by
\begin{align}\label{S2 HC transform}
    \mathcal{H}f(\lambda)\,
    =\,\int_{G}\diff{x}\,\varphi_{-\lambda}(x)\,f(x) 
    \qquad\forall\,\lambda\in\mathfrak{a},\
    \forall\,f\in\mathcal{S}(K\backslash{G/K}),
\end{align}
where $\varphi_{\lambda}\in\mathcal{C}^{\infty}(K\backslash{G/K})$ is the
spherical function of index $\lambda \in \mathfrak{a}$.
Denote by $\mathcal{S}(\mathfrak{a})^{W}$ the subspace 
of $W$-invariant functions in the Schwartz space $\mathcal{S}(\mathfrak{a})$. 
Then $\mathcal{H}$ is an isomorphism between $\mathcal{S}(K\backslash{G/K})$ 
and $\mathcal{S}(\mathfrak{a})^{W}$. The inverse spherical Fourier transform 
is given by
\begin{align}\label{S2 Inverse formula}
    f(x)\,
    =\,C_0\,\int_{\mathfrak{a}}\,\diff{\lambda}\,
        |\mathbf{c(\lambda)}|^{-2}\,
        \varphi_{\lambda}(x)\,
        \mathcal{H}f(\lambda) 
        \qquad\forall\,x\in{G},\
        \forall\,f\in\mathcal{S}(\mathfrak{a})^{W},
\end{align}
where the constant $C_0=2^{n-\ell}/(2\pi)^{\ell}|K/\mathbb{M}||W|$ depends only 
on the geometry of $\mathbb{X}$, and $|\mathbf{c(\lambda)}|^{-2}$ is the 
so-called Plancherel density. We next review some elementary facts about
the Plancherel density and the elementary spherical functions.

\subsubsection{Plancherel density}
According to the Gindikin-Karpelevič formula for the Harish-Chandra
$\mathbf{c}$-function, we can write the Plancherel density as
\begin{align*}
    |\mathbf{c}(\lambda)|^{-2}\,
    =\,\prod_{\alpha\in\Sigma_{r}^{+}}\,
        \Big|\mathbf{c}_{\alpha}
        \Big(
        \frac{\langle\alpha,\lambda\rangle}{\langle\alpha,\alpha\rangle}
        \Big)\Big|^{-2}
\end{align*}
with
\begin{align*}
    \mathbf{c}_{\alpha}(z)\,=\,
    \frac{\Gamma(\frac{\langle{\alpha,\rho}\rangle}
        {\langle{\alpha,\alpha}\rangle}
        +\frac{1}{2} m_{\alpha})}
        {\Gamma(\frac{\langle{\alpha,\rho}\rangle}
        {\langle{\alpha,\alpha}\rangle})}\,
    \frac{\Gamma(\frac{1}{2}
        \frac{\langle{\alpha,\rho}\rangle}
        {\langle{\alpha,\alpha}\rangle} 
        +\frac{1}{4} m _{\alpha} 
        + \frac{1}{2} m_{2\alpha})}
        {\Gamma(\frac{1}{2}
        \frac{\langle{\alpha,\rho}\rangle}
        {\langle {\alpha,\alpha}\rangle} 
        + \frac{1}{4} m_{\alpha})}\,
    \frac{\Gamma(iz)}
        {\Gamma(iz+ \frac{1}{2}m_{\alpha})}\,
    \frac{\Gamma(\frac{i}{2}z 
        + \frac{1}{4} m_{\alpha})}
        {\Gamma(\frac{i}{2}z 
        +\frac{1}{4} m_{\alpha} 
        + \frac{1}{2} m_{2\alpha})}.
\end{align*}
Notice that $|\mathbf{c}_{\alpha}|^{-2}$ is a differential symbol on 
$\mathbb{R}$. But the Plancherel density $|\mathbf{c}(\lambda)|^{-2}$, as a
product of one-dimensional symbols, is not a symbol on $\mathfrak{a}$ in general. 
We shall use the expressions
\begin{align*}
    |\mathbf{c}(\lambda)|^{-2}\,
    =\,\mathbf{c}(w.\lambda)^{-1}\,\mathbf{c}(-w.\lambda)^{-1}
    \qquad\forall\,\lambda\in\mathfrak{a},\,\forall\,w\in{W}
\end{align*}
as well as
\begin{align*}
    \mathbf{b}(\lambda)\,
    =\,\bm{\pi}(i\lambda)\,\mathbf{c}(\lambda)\,
    =\,\prod_{\alpha\in\Sigma_{r}^{+}}\,
        \mathbf{b}_{\alpha}
        \Big(
        \frac{\langle\alpha,\lambda\rangle}{\langle\alpha,\alpha\rangle}
        \Big)
\end{align*}
where $\bm{\pi}(i\lambda)=\prod_{\alpha\in\Sigma_{r}^{+}}
\langle{\alpha,\lambda}\rangle$ and
\begin{align*}
    \mathbf{b}_{\alpha}(z)\,
    =\,|\alpha|^{2}\,iz\,\mathbf{c}_{\alpha}(z)
    =\,|\alpha|^{2}\,
        \tfrac{\Gamma(\frac{\langle{\alpha,\rho}\rangle}
            {\langle{\alpha,\alpha}\rangle}
            +\frac{1}{2} m_{\alpha})}
            {\Gamma(\frac{\langle{\alpha,\rho}\rangle}
            {\langle{\alpha,\alpha}\rangle})}\,
        \tfrac{\Gamma(\frac{1}{2}
          \frac{\langle{\alpha,\rho}\rangle}
          {\langle{\alpha,\alpha}\rangle} 
          +\frac{1}{4} m _{\alpha} 
          + \frac{1}{2} m_{2\alpha})}
          {\Gamma(\frac{1}{2}
          \frac{\langle{\alpha,\rho}\rangle}
          {\langle {\alpha,\alpha}\rangle} 
         + \frac{1}{4} m_{\alpha})}\,
        \tfrac{\Gamma(iz+1)}
            {\Gamma(iz+ \frac{1}{2}m_{\alpha})}\,
        \tfrac{\Gamma(\frac{i}{2}z 
            + \frac{1}{4} m_{\alpha})}
            {\Gamma(\frac{i}{2}z 
            +\frac{1}{4} m_{\alpha} 
            + \frac{1}{2} m_{2\alpha})}.
\end{align*}
For every root $\alpha\in\Sigma_{r}^{+}$, the function 
$\mathbf{b}_{\alpha}(-iz)^{\pm1}$ is holomorphic for $\im{z}>-1/2$, with
\begin{align}\label{S2 behavior of Gamma functions}
    \frac{\Gamma(iz+1)}
        {\Gamma(iz+ \frac{1}{2}m_{\alpha})}\,
    \frac{\Gamma(\frac{i}{2}z 
        + \frac{1}{4} m_{\alpha})}
        {\Gamma(\frac{i}{2}z 
        +\frac{1}{4} m_{\alpha} 
        + \frac{1}{2} m_{2\alpha})}\,
    \sim\,2^{\frac{m_{2\alpha}}{2}}\,
        z^{1-\frac{m_{\alpha}}{2}-\frac{m_{2\alpha}}{2}}
    \qquad\textnormal{as}\quad\,|z|\rightarrow\infty.
\end{align}
Hence $\mathbf{b}(-\lambda)^{\pm1}$ is a holomorphic function for 
$\lambda\in\mathfrak{a}+i\overline{\mathfrak{a}^{+}}$, which has the following
behavior:
\begin{align}\label{S2 bfunction}
    |\mathbf{b}(-\lambda)|^{\pm1}\,
    \asymp\,
    \prod_{\alpha\in\Sigma_{r}^{+}}\,
        (1+|\langle{\alpha,\lambda}\rangle|)^{
        \mp\frac{m_{\alpha}+m_{2\alpha}}{2}\pm1}
\end{align}
and whose derivatives can be estimated by
\begin{align}\label{S2 bfunction derivative}
    p(\tfrac{\partial}{\partial\lambda})
    \mathbf{b}(-\lambda)^{\pm1}\,
    =\,\mathrm{O}\big(|\mathbf{b}(-\lambda)|^{\pm1}\big),
\end{align}
where $p(\tfrac{\partial}{\partial\lambda})$ is any differential 
polynomial.

\subsubsection{Spherical functions}
For every $\lambda\in\mathfrak{a}$, the spherical function
$\varphi_{\lambda}$ is a smooth bi-$K$-invariant eigenfunction of all 
$G$-invariant differential operators on $\mathbb{X}$, in particular of the
Laplace-Beltrami operator:
\begin{equation*}
    -\Delta\varphi_{\lambda}(x)\,
    =\,(|\lambda|^{2}+|\rho|^2)\,\varphi_{\lambda}(x).
\end{equation*}
It is symmetric in the sense that
$\varphi_{\lambda}(x^{-1})=\varphi_{-\lambda}(x)$,
and is given by the integral representation
\begin{align}\label{S2 Spherical Function}
    \varphi_{\lambda}(x)\, 
    =\,\int_{K}\diff{k}\,e^{\langle{i\lambda+\rho,\,A(kx)}\rangle}.
\end{align}
We list below some of the known properties of spherical functions 
$\varphi_{\lambda}$ and refer to \cite[Chap.4]{GaVa1988} and 
\cite[Chap.IV]{Hel2000} for more details.
We say that a vector $H\in\mathfrak{a}$ lies on a wall if there exists a
root $\alpha\in\Sigma$ such that $\langle{\alpha,H}\rangle=0$. 
Otherwise, we say that $H$ stays away from the walls.
For all $\lambda\in\mathfrak{a}$ \textit{regular}, 
the following converging Harish-Chandra 
expansion is known to hold for all $H\in\mathfrak{a}^{+}$ 
\begin{align}\label{S2 Harish-Chandra expansion 1}
    \varphi_{\lambda}(\exp{H})\,
    =\,\sum_{w\in{W}}\mathbf{c}(w.\lambda)\,
        \Phi_{w.\lambda}(H),
\end{align}
where
\begin{align*}
    \Phi_{\lambda}(H)\,
    =\,e^{\langle{i\lambda-\rho,H}\rangle}\,
        \sum_{q\in2Q}\gamma_{q}(\lambda)\,
        e^{-\langle{q,H}\rangle}.
\end{align*}
Here $Q=\sum_{\alpha\in\Sigma_{s}^{+}}\mathbb{N}\alpha$ denotes the positive root
lattice and the leading coefficient $\gamma_{0}$ equals $1$. For every $q\ge1$,
$\gamma_{q}(\lambda)$ are rational functions in
$\lambda\in\mathfrak{a}_{\mathbb{C}}$, which have no poles for all
$\lambda\in\mathfrak{a}+i\overline{\mathfrak{a}^{+}}$ and
satisfy there
\begin{align}\label{S2 gammaq 1}
    |\gamma_{q}(\lambda)|\,
    \lesssim\,(1+|q|)^{N_{\gamma}}
\end{align}
for some nonnegative constant $N_{\gamma}$. 
Their derivatives can be estimated as follows by using Cauchy's formula:
\begin{align}\label{S2 gammaq 2}
    p(\tfrac{\partial}{\partial\lambda})
    \gamma_{q}(\lambda)\,
    =\,\mathrm{O}\big(\gamma_{q}(\lambda)\big).
\end{align}
Moreover, all derivatives of
$\varphi_{\lambda}(\exp{H})$ in $H$ have the corresponding expansion
\begin{align}\label{S2 Harish-Chandra expansion 2}
    p(\tfrac{\partial}{\partial{H}})\,\varphi_{\lambda}(\exp{H})\,
    =\,
    e^{-\langle{\rho,H}\rangle}\,
    \sum_{q\in2Q}e^{-\langle{q,H}\rangle}\,
    \sum_{w\in{W}}\mathbf{c}(w.\lambda)\,
    \gamma_{q}(w.\lambda)\,p(iw.\lambda-\rho-q)\,
    e^{i\langle{w.\lambda,H}\rangle}.
\end{align}
Recall that all the elementary spherical functions $\varphi_{\lambda}$ 
with parameter $\lambda\in\mathfrak{a}$ are controlled by the ground spherical
function $\varphi_{0}$, which satisfies the global estimate
\begin{align}\label{S2 global estimate phi0}
    \varphi_{0}(\exp{H})\,
    \asymp\,
        \Big\lbrace \prod_{\alpha\in\Sigma_{r}^{+}} 
        1+\langle\alpha,H\rangle\Big\rbrace\,
        e^{-\langle\rho, H\rangle}
    \qquad\forall\,H\in\overline{\mathfrak{a}^{+}}.
\end{align}

\subsection{Heat kernel on symmetric spaces}
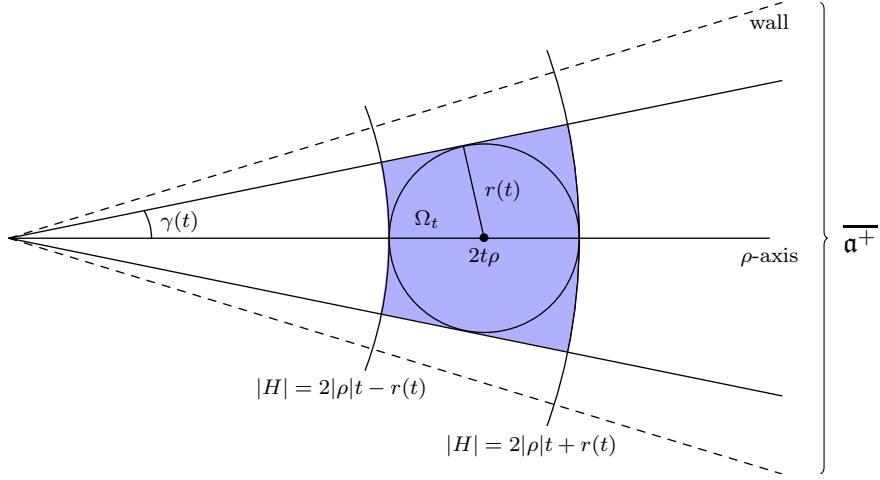
\begin{figure}
    \centering
    \input{Concentration1.tex}
    \caption{Critical region where the heat kernel concentrates.}
    \label{fig concentration1}
\end{figure}

The heat kernel on $\mathbb{X}$ is a positive bi-$K$-invariant right 
convolution kernel, i.e., $h_{t}(xK,yK)=h_{t}(y^{-1}x)>0$, 
which is thus determined by its restriction 
to the positive Weyl chamber. 
According to the inversion formula of the spherical Fourier transform, 
the heat kernel is given by
\begin{align}\label{S2 heat kernel inv}
    h_{t}(x)\,
    =\,C_{0}\,\int_{\mathfrak{a}}\,\diff{\lambda}\,
        |\mathbf{c(\lambda)}|^{-2}\,
        \varphi_{\lambda}(x)\,
        e^{-t(|\lambda|^{2}+|\rho|^{2})}
\end{align}
and satisfies the global estimate
\begin{align}\label{S2 heat kernel}
    h_{t}(\exp{H})\,
    \asymp\,t^{-\frac{n}{2}}\,
    \Big\lbrace{
    \prod_{\alpha\in\Sigma_{r}^{+}}
    (1+t+\langle{\alpha,H}\rangle)^{\frac{m_{\alpha}+m_{2\alpha}}{2}-1}
    }\Big\rbrace\,\varphi_{0}(\exp{H})
    e^{-|\rho|^{2}t-\frac{|H|^{2}}{4t}}
\end{align}
for all $t>0$ and $H\in\overline{\mathfrak{a}^{+}}$, 
see \cite{AnJi1999,AnOs2003}.
Recall that $\int_{\mathbb{X}}\diff{x}\,h_{t}(x)=1$.
Let $r(t)=\frac{\sqrt{t}}{\varepsilon(t)}$, 
where $\varepsilon(t)$ is a positive function, such that 
$\varepsilon(t)\searrow0$ and $\sqrt{t}\,\varepsilon(t)\rightarrow\infty$,
as $t\rightarrow\infty$.
Let $\gamma(t)$ be another positive function such that 
$\sqrt{t}\gamma(t)\rightarrow\infty$ as $t\rightarrow\infty$.
Denote by $\Omega_{\textrm{annulus}}\in\mathfrak{a}$ the annulus 
$2|\rho|t-r(t)\le|H|\le2|\rho|t+r(t)$ and by 
$\Omega_{\textrm{cone}}\in\mathfrak{a}$ the solid cone with angle
$\gamma(t)$ around the $\rho$-axis. 
According to \cite{AnSe1992}, the heat kernel $h_{t}$ is asymptotically 
concentrated in $K(\exp(\Omega_{\textrm{annulus}}\cap\Omega_{\textrm{cone}}))K$
(see the blue region in \cref{fig concentration1}):
\begin{align*}
    \int_{K(\exp(\Omega_{\textrm{annulus}}\,\cap\,\Omega_{\textrm{cone}}))K}
    \diff{x}\,h_{t}(x)\,
    \longrightarrow\,1
    \qquad\textnormal{as}\quad\,t\rightarrow\infty.
\end{align*}
The following lemma shows that the heat kernel 
$h_{t}(x)$ concentrates indeed in $K(\exp\Omega_{t})K$, 
where $\Omega_{t}=B(2t\rho,r(t))$ 
denotes the ball in $\mathfrak{a}$ with center $2t\rho$
and radius $r(t)$.

\begin{lemma}\label{S2 prop concentration}
The following estimate holds for any $N\ge0$:
\begin{align}\label{S2 kernel estimate outside}
    \|h_{t}\|_{L^{1}(G\smallsetminus{K(\exp\Omega_{t})K})}\,
    \lesssim\,\varepsilon(t)^{N}
    \qquad\forall\,t>1.
\end{align}
\end{lemma}

\begin{proof}
For simplicity, we denote by $H=x^{+}$ the middle component of $x\in{G}$ in the 
Cartan decomposition. Let $\eta>0$ be a constant such that $\eta{t}>{r(t)}$.
The proof of \eqref{S2 kernel estimate outside}
is based on the following large time heat kernel estimates, 
which follow from \eqref{S2 heat kernel} for some $L\ge0$:
\begin{align}
\begin{split}
    h_{t}(\exp{H})\,\lesssim\,
        \underbrace{\vphantom{\Big|}
        e^{-|\rho|^{2}t-\langle{\rho,H}\rangle-\frac{|H|^{2}}{4t}}
        }_{e^{-\frac{|H-2t\rho|^{2}}{4t}-2\langle{\rho,H}\rangle}}
    \times\,
    \begin{cases}
        t^{-\frac{\ell}{2}}
        &\qquad\textrm{if}\,\,\,H\in{B(2t\rho,\eta{t})},\\[5pt]
        t^{L}
        &\qquad\textrm{if}\,\,\,
        H\in{B(0,3|\rho|t})\smallsetminus{B(2t\rho,\eta{t})},\\[5pt]
        |H|^{L}
        &\qquad\textrm{if}\,\,\,H\ge3|\rho|t.
    \end{cases}
\end{split} 
\end{align}
Firstly, we have
\begin{align*}
    \int_{K(\exp{B(2t\rho,\eta{t})})K\,
        \smallsetminus\,K(\exp\Omega_{t})K}
    \diff{x}\,h_{t}(x)\,
    &\lesssim\,\int_{r(t)\,\le\,|H-2\rho{t}|\,\le\,\eta{t}}\diff{H}\,
        \delta(H)\,h_{t}(\exp{H})\\[5pt]
    &\lesssim\,t^{-\frac{\ell}{2}}\,
        \int_{r(t)\,\le\,|H-2\rho{t}|\,\le\,\eta{t}}\diff{H}\,
        e^{-\frac{|H-2t\rho|^{2}}{4t}}\\[5pt]
    &\lesssim\,t^{-\frac{\ell}{2}}\,
        \int_{r(t)}^{\eta{t}}\diff{R}\,R^{\ell-1}\,
        e^{-\frac{R^{2}}{4t}}\\[5pt]
    &=\,\int_{\frac{r(t)}{2\sqrt{t}}}^{\frac{\eta\sqrt{t}}{2}}\diff{R}\,
        R^{\ell-1}\,e^{-R^2}\,
        \lesssim\,\big(\tfrac{\sqrt{t}}{r(t)}\big)^{N}
\end{align*}
for any $N\ge0$.
Similar computations imply secondly 
\begin{align*}
    \int_{K(\exp{B(0,3|\rho|t)})K\,
        \smallsetminus\,K(\exp{B(2t\rho,\eta{t})})K}
        \diff{x}\,h_{t}(x)\,
    &\lesssim\,
        \int_{[\,\overline{\mathfrak{a}^{+}}\cap{B(0,3|\rho|t)}]\,
        \smallsetminus\,B(2t\rho,\eta{t})}\diff{H}\,
        \delta(H)\,h_{t}(\exp{H})\\[5pt]
    &\lesssim\,t^{L}\,
        \int_{\eta{t}\,\le\,|H-2\rho{t}|\,\le\,5|\rho|t}\diff{H}\,
        e^{-\frac{|H-2t\rho|^{2}}{4t}}\\[5pt]
    &\lesssim\,t^{L}\,
        \int_{\eta{t}}^{5|\rho|t}\diff{R}\,R^{\ell-1}\,
        e^{-\frac{R^{2}}{4t}}\\[5pt]
    &=\,t^{L+\frac{\ell}{2}}\,
        \int_{\frac{\eta\sqrt{t}}{2}}^{\frac{5|\rho|\sqrt{t}}{2}}
        \diff{R}\,R^{\ell-1}\,e^{-R^2}\,
        \lesssim\,t^{-N}
\end{align*}
and thirdly
\begin{align*}
    \int_{G\,\smallsetminus\,K(\exp{B(0,3|\rho|t)})K}
        \diff{x}\,h_{t}(x)\,
    &\lesssim\,
        \int_{\,\overline{\mathfrak{a}^{+}}\,
        \smallsetminus\,B(0,3|\rho|t)}\diff{H}\,
        \delta(H)\,h_{t}(\exp{H})\\[5pt]
    &\lesssim\,
        \int_{|H|\,\ge\,3|\rho|t}\diff{H}\,|H|^{L}\,
        e^{-\frac{|H-2t\rho|^{2}}{4t}}\\[5pt]
    &\lesssim\,
        \int_{|H-2t\rho|\,\ge\,|\rho|t}\diff{H}\,|H-2t\rho|^{L}\,
        e^{-\frac{|H-2t\rho|^{2}}{4t}}\\[5pt]
    &\lesssim\,
        \int_{|\rho|t}^{+\infty}\diff{R}\,R^{L+\ell-1}\,
        e^{-\frac{R^{2}}{4t}}\,
        \lesssim\,t^{-N}
\end{align*}
for any $N\ge0$. The proposition follows by combining these three estimates.
\end{proof}

As the vectors in $\Omega_t$ stay far away from walls, the large time
behavior of the heat kernel can be described more accurately by the following
asymptotics \cite[Theorem 5.1.1]{AnJi1999}:
\begin{align}\label{S2 heat kernel critical region}
    h_t(\exp{H})\,
    \sim\,C_{1}\,t^{-\frac{\nu}{2}}\,
        \mathbf{b}\big(-i\tfrac{H}{2t}\big)^{-1}\,
        \varphi_{0}(\exp{H})\,e^{-|\rho|^{2}t-\frac{|H|^{2}}{4t}}
\end{align}
as $t\rightarrow\infty$ and 
$\mu(H)=\min_{\alpha\in\Sigma^{+}}\langle{\alpha,H}\rangle\rightarrow\infty$.
Here $C_{1}=C_{0}2^{-|\Sigma_{r}^{+}|}|W|\pi^{\frac{\ell}{2}}
\bm{\pi}(\rho_{0})\mathbf{b}(0)^{-1}$ is a positive constant, and $\rho_{0}$
denotes the half sum of all positive reduced roots.
Moreover, the ground spherical function satisfies 
\begin{align}\label{S2 phi0 far}
    \varphi_{0}(\exp{H})\,
    \sim\,C_{2}\,\bm{\pi}(H)\,e^{-\langle{\rho,H}\rangle}
\end{align}
as $\mu(H)\rightarrow\infty$, where $C_{2}=\bm{\pi}(\rho_{0})^{-1}\mathbf{b}(0)$,
see for instance \cite[Proposition 2.2.12.(ii)]{AnJi1999}.

\section{Asymptotic convergence associated with the Laplace-Beltrami operator
}\label{Section.3 Asymp}
Consider the heat equation
\begin{align}\label{HE X}
    \partial_{t}u(t,x)\,=\,\Delta_{x}u(t,x),
    \qquad
    u(0,x)\,=\,u_{0}(x),
\end{align}
with $u_{0}\in{L}^{1}(\mathbb{X})$
and denote by $M$ the mass of the initial data:
\begin{align*}
    M\,
    =\,\int_{G}\diff{x}\,u_{0}(x).
\end{align*}
In this section, we prove our first main result \cref{S1 Main thm 1}.
We start with smooth bi-$K$-invariant compactly supported initial data and study 
the long-time asymptotic convergence in the critical region 
where the heat kernel concentrates. 
By completing the estimates outside the critical region and by using a 
standard density argument, we next establish \eqref{S1 Main thm 1 convergence}
for all bi-$K$-invariant $L^1$ initial data.
Finally, we provide a counterexample in the non bi-$K$-invariant case.

\subsection{Heat asymptotics in the critical region for 
$\mathcal{C}_{c}^{\infty}(K\backslash{G/K})$ initial data}\label{S3 Sub1}

In order to treat heat asymptotics in the critical region, we pass to the frequency side. For an alternative proof of the following \cref{S3 L1 convergence in critical region},  directly on the space side, see \cref{Remark 3.10}.

As both the heat kernel $h_{t}$ and the initial data $u_0$ are bi-$K$-invariant,
we have $\mathcal{H}(u_0\ast h_t)=(\mathcal{H}u_0)(\mathcal{H}h_t)$, 
see for instance \cite[p. 347]{Ank1991}. 
Then, we can write the solution to \eqref{HE X} as
\begin{align*}
    u(t,x)\,=\,(u_{0}*h_{t})(x)
    =\,C_{0}\,\int_{\mathfrak{a}}\diff{\lambda}\,
        |\mathbf{c}(\lambda)|^{-2}\,
        e^{-t(|\lambda|^2+|\rho|^{2})}\,\varphi_{\lambda}(x)\,
        \mathcal{H}u_{0}(\lambda)
\end{align*}
by using the inversion formula \eqref{S2 Inverse formula} 
of the spherical Fourier transform and \eqref{S2 heat kernel inv}.
Notice that the mass $M$ equals $\mathcal{H}u_{0}(-i\rho)$ 
by the spherical Fourier transform and \eqref{S2 Spherical Function}. 
Our aim in this subsection is to estimate the difference
\begin{align}\label{S3 u-Mht}
    u(t,x)\,-Mh_{t}(x)\,
    =\,C_{0}\,e^{-|\rho|^2t}\,
       \int_{\mathfrak{a}}\diff{\lambda}\,
       |\mathbf{c}(\lambda)|^{-2}\,
       e^{-t|\lambda|^2}\,\varphi_{\lambda}(x)\,
       \underbrace{\vphantom{\Big|}
       \lbrace{
       \mathcal{H}u_{0}(\lambda)-\mathcal{H}u_{0}(-i\rho)
       }\rbrace
       }_{U(\lambda)}
\end{align}
in the critical region $K(\exp\Omega_{t})K$.
Recall that $\Omega_{t}=B(2t\rho,r(t))$ where $r(t)$ grows slightly faster 
than $\sqrt{t}$ and slower than $t$ as $t\rightarrow\infty$. 
To this end we resume the computations in \cite[Step 2, pp. 1054-1056]{AnJi1999}.

\begin{proposition}\label{S3 L1 convergence in critical region}
Let $u_{0}$ be a smooth bi-$K$-invariant initial data with compact support. Then
\begin{align}\label{S3 convergene critical region}
    \|u(t,\,\cdot\,)-Mh_{t}\|_{L^{1}(K(\exp\Omega_{t})K)}\,
    \longrightarrow\,0
    \qquad\textrm{as}\,\,\,t\rightarrow\,0.
\end{align}
\end{proposition}

\begin{proof}
For simplicity, we denote by $H=x^{+}$ the middle component of $x\in{G}$ in the 
Cartan decomposition. By substituting \eqref{S2 Harish-Chandra expansion 1} 
in \eqref{S3 u-Mht} and by writing $|\mathbf{c}(\lambda)|^{-2}
=\mathbf{c}(\lambda)^{-1}\mathbf{c}(-\lambda)^{-1}$, 
with $\mathbf{c}(-\lambda)^{-1}=\bm{\pi}(-i\lambda)\mathbf{b}(-\lambda)^{-1}$, 
we obtain first
\begin{align}\label{S3 u-Mht H}
    u(t,\exp{H})\,-Mh_{t}(\exp{H})\,
    =\,C_{0}\,|W|\,e^{-|\rho|^2t-\langle{\rho,H}\rangle}\,
        \sum_{q\in{Q}}e^{-\langle{q,H}\rangle}\,E_{q}(t,H)
\end{align}
where
\begin{align*}
    E_{q}(t,H)\,
    =\,\int_{\mathfrak{a}}\diff{\lambda}\,
        \bm{\pi}(-i\lambda)\,e^{-t|\lambda|^{2}}\,
        e^{i\langle{\lambda,H}\rangle}\,\mathbf{b}(-\lambda)^{-1}\,
        \gamma_{q}(\lambda)\,U(\lambda).
\end{align*}
By using the formula
\begin{align*}
    \bm{\pi}(-i\lambda)\,e^{-t|\lambda|^{2}}\,
    =\,\bm{\pi}(\tfrac{i}{2t}\tfrac{\partial}{\partial\lambda})\,
        e^{-t|\lambda|^{2}}
\end{align*}
and by integrating by parts, we obtain firstly
\begin{align*}
    E_{q}(t,H)\,
    &=\,(2t)^{-|\Sigma_{r}^{+}|}\,
        \sum_{\Sigma_{r}^{+}=\Sigma'\sqcup\Sigma''}
        \Big\lbrace{
        \prod_{\alpha'\in\Sigma'}\langle{\alpha',H}\rangle
        }\Big\rbrace\\[5pt]
    &\times\,\int_{\mathfrak{a}}\diff{\lambda}\,
        e^{-t|\lambda|^{2}}\,e^{i\langle{\lambda,H}\rangle}\,
        \Big\lbrace{
        \prod_{\alpha''\in\Sigma''}(-i\partial_{\alpha''})
        }\Big\rbrace
        \big\lbrace{
        \mathbf{b}(-\lambda)^{-1}\,\gamma_{q}(\lambda)\,U(\lambda)
        }\big\rbrace.
\end{align*}
After moving the contour of integration according to 
$\lambda\mapsto\lambda+i\frac{H}{2t}$ and rescaling the integral according to
$\lambda\mapsto\frac{\lambda}{\sqrt{t}}$, we obtain secondly
\begin{align}\label{S3 Eq}
\begin{split}
    E_{q}(t,H)\,
    &=\,2^{-|\Sigma_{r}^{+}|}\,t^{-\frac{\nu}{2}}\,e^{-\frac{|H|^{2}}{4t}}
        \sum_{\Sigma_{r}^{+}=\Sigma'\sqcup\Sigma''}
        \Big\lbrace{
        \prod_{\alpha'\in\Sigma'}\langle{\alpha',H}\rangle
        }\Big\rbrace\\[5pt]
    &\times\,\int_{\mathfrak{a}}\diff{\lambda}\,e^{-|\lambda|^{2}}\,
        \Big\lbrace{
        \prod_{\alpha''\in\Sigma''}(-i\sqrt{t}\partial_{\alpha''})
        }\Big\rbrace
        \big\lbrace{
        \mathbf{b}\big(\!-\tfrac{\lambda}{\sqrt{t}}-i\tfrac{H}{2t}\big)^{-1}\,
        \gamma_{q}\big(\tfrac{\lambda}{\sqrt{t}}+i\tfrac{H}{2t}\big)\,
        U\big(\tfrac{\lambda}{\sqrt{t}}+i\tfrac{H}{2t}\big)
        }\big\rbrace.
\end{split}
\end{align}
Now let us study \eqref{S3 Eq} in different cases.
We first assume that no derivative in \eqref{S3 Eq} hits $U$. 
According to the Paley-Wiener theorem 
(see for instance \cite[Chap.III, Theorem 5.1]{Hel1994}), 
as $u_{0}$ is a smooth bi-$K$-invariant
function with compact support on $G$, $\mathcal{H}u_{0}$ is a $W$-invariant 
holomorphic function on $\mathfrak{a}_{\mathbb{C}}$ such that
\begin{align}\label{S3 Paley-Wiener}
    \exists\,C'\ge0,\,
    \forall\,j\in\mathbb{N},\,\forall\,N\in\mathbb{N},\,
    \exists\,C''\ge0,\,
    \forall\,\lambda\in\mathfrak{a}_{\mathbb{C}},\,
    |\nabla_{\lambda}^{j}\mathcal{H}u_{0}(\lambda)|\,
    \le\,C''(1+|\lambda|)^{-N}e^{C'|\im\lambda|}.
\end{align}
Hence
\begin{align}\label{S3 estim U}
    |U\big(\tfrac{\lambda}{\sqrt{t}}+i\tfrac{H}{2t}\big)|\,
    \lesssim\,e^{\frac{C'}{2}\frac{|H|}{t}}\,
        \big(\tfrac{|\lambda|}{\sqrt{t}}+|\tfrac{H}{2t}-\rho|\big)
    \qquad\forall\,\lambda,\,H\in\mathfrak{a},
    \,\forall\,t>0.
\end{align}
On the other hand, according to \eqref{S2 bfunction} and  
\eqref{S2 bfunction derivative},
\begin{align}\label{S3 estim b}
    \big(\sqrt{t}\,
    \nabla_{\lambda}\big)^{j}\,
    \mathbf{b}\big(\!-\tfrac{\lambda}{\sqrt{t}}-i\tfrac{H}{2t}\big)^{-1}\,
    =\,\mathrm{O}\Big(
        \big(1+\tfrac{|\lambda|}{\sqrt{t}}\big)^{m}\,
        \prod\nolimits_{\alpha\in\Sigma_{r}^{+}}\,
        (1+\tfrac{\langle{\alpha,H}\rangle}{t})^{
        (\frac{m_{\alpha}+m_{2\alpha}}{2}-1)_{+}}\Big)
\end{align}
where $m=\sum_{\alpha\in\Sigma_{r}^{+}}(\frac{m_{\alpha}+m_{2\alpha}}{2}-1)_{+}$,
and according to \eqref{S2 gammaq 1} and  \eqref{S2 gammaq 2},
\begin{align}\label{S3 estim gammaq}
    \big(\sqrt{t}\tfrac{\partial}{\partial\lambda}\big)^{j}\,
    \gamma_{q}\big(\tfrac{\lambda}{\sqrt{t}}+i\tfrac{H}{2t}\big)\,
    =\,\mathrm{O}\big((1+|q|)^{N_{\gamma}}\big),
\end{align}
where $N_{\gamma}$ is a positive constant.
By using \eqref{S3 estim U}, \eqref{S3 estim b} and \eqref{S3 estim gammaq},
we estimate the terms in \eqref{S3 Eq} where no derivative hits $U$ by
\begin{align}\label{S3 estim no derivative hits U}
    t^{-\frac{\nu}{2}}\,e^{-\frac{|H|^{2}}{4t}}\,
    \overbrace{\vphantom{\Big|}
        (1+|H|)^{|\Sigma_{r}^{+}|}}^{\asymp\,t^{|\Sigma_{r}^{+}|}}\,
    &\Big\lbrace{\prod\nolimits_{\alpha\in\Sigma_{r}^{+}}\,
        \overbrace{\vphantom{\Big|}
            (1+\tfrac{\langle{\alpha,H}\rangle}{t})^{
            (\frac{m_{\alpha}+m_{2\alpha}}{2}-1)_{+}}}^{\asymp\,1}
    }\Big\rbrace\,
    (1+|q|)^{N_{\gamma}}\,
    \overbrace{\vphantom{\Big|}
        e^{\frac{C'}{2}\frac{|H|}{t}}}^{\lesssim\,1}\,
    \overbrace{\vphantom{\Big|}
    \big(\tfrac{1}{\sqrt{t}}+\tfrac{r(t)}{t}\big)}^{\lesssim\,\tfrac{r(t)}{t}}
    \nonumber\\[5pt]
    &\lesssim\,(1+|q|)^{N_{\gamma}}\,
        t^{-\frac{\ell}{2}-1}\,r(t)\,e^{-\frac{|H|^{2}}{4t}}
\end{align}
when $t$ is large and $H\in\Omega_{t}$.
If some derivatives $\partial_{\alpha''}$ in \eqref{S3 Eq} hit $U$, then
\begin{align*}
    \Big|\prod_{\alpha'\in\Sigma'}\langle{\alpha',H}\rangle\Big|\,
    \lesssim\,|H|^{|\Sigma'|}\,
    \lesssim\,t^{|\Sigma_{r}^{+}|-1}
\end{align*}
and, by using \eqref{S3 Paley-Wiener}, \eqref{S3 estim b} and 
\eqref{S3 estim gammaq}, we obtain this time that \eqref{S3 Eq} is bounded by
\begin{align}\label{S3 estim a derivative hits U}
    (1+|q|)^{N_{\gamma}}\,t^{-\frac{\ell}{2}-1}\,e^{-\frac{|H|^{2}}{4t}}.
\end{align}
By combining \eqref{S3 u-Mht H}, \eqref{S3 Eq}, 
\eqref{S3 estim no derivative hits U} and \eqref{S3 estim a derivative hits U},
we obtain
\begin{align*}
    |u(t,\exp{H})\,-Mh_{t}(\exp{H})|\,
    \lesssim\,
    \underbrace{\vphantom{\Big|}\Big\lbrace{
    \sum_{q\in2Q}(1+|q|)^{N_{\gamma}}\,e^{-\langle{q,H}\rangle}
    \Big\rbrace}}_{\lesssim\,1}
    t^{-\frac{\ell}{2}-1}\,r(t)\,
    e^{-|\rho|^{2}t-\langle{\rho,H}\rangle-\frac{|H|^{2}}{4t}}
\end{align*}
when $t$ is large and $H\in\Omega_{t}$. 
By integration, we conclude that
\begin{align*}
    \int_{K(\exp\Omega_{t})K}\diff{x}\,|u(t,x)-Mh_{t}(x)|
    &=\,C_{0}\,\int_{\Omega_{t}}\diff{H}\,
        \delta(H)\,|u(t,\exp{H})-Mh_{t}(\exp{H})|\\[5pt]
    &\lesssim\,t^{-\frac{\ell}{2}-1}\,r(t)\,
        \int_{\Omega_{t}}\diff{H}\,
        e^{-|\rho|^{2}t+\langle{\rho,H}\rangle-\frac{|H|^{2}}{4t}}
        \\[5pt]
    &=\,\tfrac{r(t)}{t}\,
        t^{-\frac{\ell}{2}}\,\int_{B(2t\rho,r(t))}\diff{H}\,
        e^{-\frac{|H-2t\rho|^{2}}{4t}}\\[5pt]
    &=\,\tfrac{r(t)}{t}\,
        \underbrace{\vphantom{\Big|}
        \int_{0}^{\frac{r(t)}{\sqrt{t}}}\diff{r}\,
        r^{\ell-1}\,e^{-\frac{r^{2}}{4}}
        }_{<\infty}
        \lesssim\,\tfrac{r(t)}{t}
\end{align*}
tends to $0$ at speed $\frac{r(t)}{t}=\tfrac{1}{\sqrt{t}\,\varepsilon(t)}$.
\end{proof}       
\subsection{Estimates outside the critical region}\label{S3 Sub2}
In this subsection, we show that the solution $u(t,x)$ to the heat equation 
\eqref{HE X} vanishes asymptotically 
in $L^{1}(G\smallsetminus{K(\exp\Omega_{t})K})$ as $t\rightarrow\infty$.

\begin{proposition}
The solution to the heat equation \eqref{HE X} satisfies
\begin{align}\label{S3 solution estimate outside}
    \|u(t,\,\cdot\,)\|_{L^{1}(G\smallsetminus{K(\exp\Omega_{t})K})}\,
    \lesssim\,\varepsilon(t)^{N}
\end{align}
for $t>0$ large enough and for any $N\ge0$.
\end{proposition}

\begin{lemma}\label{lemma distance}
$d(K(\exp{H_{1}})K,K(\exp{H_{2}})K)=|H_{1}-H_{2}|$
for all $H_{1},H_{2}\in\overline{\mathfrak{a}^{+}}$.
\end{lemma}

\begin{proof}
On the one hand,
\begin{align*}
    d(K(\exp{H_{1}})K,K(\exp{H_{2}})K)\,
    =\,\inf_{k_{1},k_{2}\in{K}}d(k_{1}(\exp{H_{1}})K,k_{2}(\exp{H_{2}})K)
    \ge\,|H_{1}-H_{2}|
\end{align*}
according to \eqref{S2 Distance}. On the other hand,
\begin{align*}
    d(K(\exp{H_{1}})K,K(\exp{H_{2}})K)\,
    &\le\,d((\exp{H_{1}})K,(\exp{H_{2}})K)\\[5pt]
    &=\,d(\exp{(H_{1}-H_{2}})K,eK)\,
    =\,|H_{1}-H_{2}|
\end{align*}
as $(\exp\mathfrak{a})K$ is a flat subspace of $G/K$.
\end{proof}

\begin{proof}[Proof of Proposition 3.2]
Since $h_{t}$ and $u_{0}$ are both bi-$K$-invariant functions on $G$,
we have
\begin{align*}
    u(t,x)\,
    =\,(u_{0}*h_{t})(x)\,
    =\,(h_{t}*u_{0})(x)\,
    =\,\int_{G}\diff{y}\,h_{t}(xy^{-1})\,u_{0}(y).
\end{align*}
Let $\xi>0$ be a constant such that the compact support of $u_{0}$
belongs to $K(\exp{B(0,\xi)})K$. Then
\begin{align*}
    \int_{G\,\smallsetminus\,K(\exp\Omega_{t})K}\diff{x}\,|u(t,x)|\,
    \lesssim\,
    \int_{K(\exp{B(0,\xi)})K}\diff{y}\,|u_{0}(y)|\,
    \int_{G\,\smallsetminus\,K(\exp\Omega_{t})K}\diff{x}\,h_{t}(xy^{-1}).
\end{align*}
We notice from the \cref{lemma distance} that
\begin{align*}
    K(\exp\Omega_{t})K\,
    =\,K(\exp{B(2t\rho,r(t)}))K\,
    =\,\lbrace{
    zK\in{G/K}\,|\,d(zK,\,K(\exp{2t\rho})K)<r(t)
    }\rbrace.
\end{align*}
Hence
\begin{align*}
    x\,\in\,G\smallsetminus{K(\exp\Omega_{t})K}
    \qquad\Longrightarrow\qquad
    xy^{-1}\,\in\,G\smallsetminus{K(\exp{B(2t\rho,r(t)-\xi}))K}.
\end{align*}
Indeed, if $xy^{-1}\in{K(\exp{B(2t\rho,r(t)-\xi}))K}$, then
\begin{align*}
    d(xK,\,K(\exp2t\rho)K)\,
    \le\,\underbrace{\vphantom{\Big|}
        d(xK,\,xy^{-1}K)}_{d(yK,\,eK)\,<\,\xi}\,
    +\,\underbrace{\vphantom{\Big|}
        d(xy^{-1}K,\,K(\exp2t\rho)K)}_{<\,r(t)-\xi}\,
    <\,r(t),
\end{align*}
which implies that $x\in{K(\exp\Omega_{t})K}$. Therefore
\begin{align}\label{S3 htz}
    \int_{G\,\smallsetminus\,K(\exp\Omega_{t})K}\diff{x}\,h_{t}(xy^{-1})\,
    \le\,
    \int_{G\,\smallsetminus\,K(\exp{B(2t\rho,r(t)-\xi}))K}\diff{z}\,
        h_{t}(z)\,
\end{align}
for all $y\in{K(\exp{B(0,\xi)})K}$. 
As $r(t)-\xi\ge\frac{r(t)}{2}$ for $t$ large enough, we deduce from 
\eqref{S2 kernel estimate outside} that the right-hand side of \eqref{S3 htz}
is also $\textrm{O}(\varepsilon(t)^{N})$. In conclusion,
\begin{align*}
    \int_{G\,\smallsetminus\,K(\exp\Omega_{t})K}\diff{x}\,|u(t,x)|\,
    \lesssim\,
    \varepsilon(t)^{N}
    \qquad\forall\,N\ge0.
\end{align*}
\end{proof}

\begin{remark}
The estimate \eqref{S3 solution estimate outside} still holds without the
bi-$K$-invariance assumption. Indeed, by using the symmetries 
$(K(\exp\Omega_{t})K)^{-1}=K(\exp\Omega_{t})K$ and 
$h_{t}(y^{-1}x^{-1})=h_{t}(xy)$, we have
\begin{align*}
    \int_{G\,\smallsetminus\,K(\exp\Omega_{t})K}\diff{x}\,|u(t,x)|\,
    &=\,
    \int_{G\,\smallsetminus\,K(\exp\Omega_{t})K}\diff{x}\,|u(t,x^{-1})|\\[5pt]
    &\lesssim\,
    \int_{K(\exp{B(0,\xi)})K}\diff{y}\,|u_{0}(y)|\,
    \int_{G\,\smallsetminus\,K(\exp\Omega_{t})K}\diff{x}\,h_{t}(xy)
\end{align*}
and we can conclude similarly.
\end{remark}

\subsection{Long-time convergence for general bi-$K$-invariant data}
In the previous two subsections, we studied the long-lime asymptotic behavior 
for the heat equation with smooth compactly supported bi-$K$-invariant initial
data. 
Using those estimates and a standard density argument, 
we first prove in this subsection \cref{S1 Main thm 1} for the whole class of
$L^{1}(\mathbb{X}$) functions that are bi-$K$-invariant.
Then we discuss the same problem in $L^{p}(\mathbb{X})$ for $p>1$.

\begin{proof}[Proof of \cref{S1 Main thm 1}]
Let $\varepsilon>0$ and $U_{0}\in\mathcal{C}_{c}^{\infty}(K\backslash{G}/K)$ 
be such that $\|u_{0}-U_{0}\|_{L^{1}(\mathbb{X})}<\tfrac{\varepsilon}{3}$.
Denote by $M_{U}=\int_{\mathbb{X}}\diff{x}\,U_{0}(x)$ the 
mass of $U_{0}$, then
\begin{align*}
    |M-M_{U}|\,
    \le\,\|u_{0}-U_{0}\|_{L^{1}(\mathbb{X})}\,
    <\,\tfrac{\varepsilon}{3}.
\end{align*}
Let $U(t,x)=(U_{0}*h_{t})(x)$ be the solution to
the heat equation with initial data $U_{0}$. We deduce from
\eqref{S2 kernel estimate outside}, \eqref{S3 convergene critical region} and
\eqref{S3 solution estimate outside} that, there exists $T>0$ such that
\begin{align*}
    \|U(t,\,\cdot\,)-M_{U}h_{t}\|_{L^{1}(\mathbb{X})}
    &\le\,\|U(t,\,\cdot\,)-M_{U}h_{t}\|_{L^{1}(K(\exp\Omega_{t})K)}\,
    +\,\|U(t,\,\cdot\,)\|_{L^{1}(G\smallsetminus{K(\exp\Omega_{t})K})}\\[5pt]
    &+\,|M_{U}|\,
        \|h_{t}\|_{L^{1}(G\smallsetminus{K(\exp\Omega_{t})K})}\\[5pt]
    &<\,\tfrac{\varepsilon}{3}
\end{align*}
for all $t\ge{T}$. In conclusion,
\begin{align*}
    \|u(t,\,\cdot\,)-Mh_{t}\|_{L^{1}(\mathbb{X})}
    &\le\overbrace{\vphantom{\Big|}
    \|u(t,\,\cdot\,)-U(t,\,\cdot\,)\|_{L^{1}(\mathbb{X})}
    }^{\le\,\|u_{0}-U_{0}\|_{L^{1}}\,\|h_{t}\|_{L^{1}}}
    +\overbrace{\vphantom{\Big|}
    \|M_{U}h_{t}-Mh_{t}\|_{L^{1}(\mathbb{X})}
    }^{\le\,|M-M_{U}|\,\|h_{t}\|_{L^{1}}}\\[5pt]
    &+\|U(t,\,\cdot\,)-M_{U}h_{t}\|_{L^{1}(\mathbb{X})}
    \\[5pt]
    &<\,\tfrac{\varepsilon}{3}+\tfrac{\varepsilon}{3}+\tfrac{\varepsilon}{3}\,
    =\varepsilon
\end{align*}
for all $\varepsilon>0$ and $t$ large enough.
\end{proof}

\vspace{5pt}
Let us turn to the long-time convergence in $L^{p}(\mathbb{X})$ with $p>1$. We first deal with the case $p=\infty$ and reach the full range by convexity.
The kernel estimate \eqref{S2 heat kernel} gives us the sup norm estimate:
\begin{align*}
    \|u(t,\,\cdot\,)-Mh_{t}\|_{L^{\infty}(\mathbb{X})}\,
    \le\,
    \|u_{0}\|_{L^{1}(\mathbb{X})}\,\|h_{t}\|_{L^{\infty}(\mathbb{X})}\,
    +\,
    |M|\,\|h_{t}\|_{L^{\infty}(\mathbb{X})}\,
    \lesssim\,t^{-\frac{\nu}{2}}e^{-|\rho|^{2}t}
\end{align*}
for $t$ large and for all $u_{0}\in{L^{1}(\mathbb{X})}$,
see for instance \cite[Proposition 4.1.1]{AnJi1999}.
In other words, we have
\begin{align}\label{S3 Linfty convergence}
    \|u(t,\,\cdot\,)-Mh_{t}\|_{L^{\infty}(\mathbb{X})}\,
    =\,
    \mathrm{O}\big(t^{-\frac{\nu}{2}}e^{-|\rho|^{2}t}\big)
    \qquad\textnormal{as}\quad\,t\rightarrow\infty.
\end{align}
Notice that such an estimate holds without the bi-$K$-invariant assumption.
By convexity, we obtain the following estimates in the $L^{p}(\mathbb{X})$
setting.
\begin{corollary}
Under the assumption of \cref{S1 Main thm 1}, we have
\begin{align}\label{S3 Lp convergence}
    \|u(t,\,\cdot\,)-Mh_{t}\|_{L^{p}(\mathbb{X})}\,
    =\,\mathrm{o}\big(t^{-\frac{\nu}{2p'}}
        e^{-\frac{|\rho|^{2}t}{p'}}\big)
    \qquad\textnormal{as}\quad\,t\rightarrow\infty
\end{align}
for all $1<p<\infty$.
\end{corollary}

\begin{remark}
As previously observed by Vázquez on real hyperbolic spaces, 
the $L^p$ norm estimate \eqref{S3 Lp convergence} is optimal, 
but the sup norm estimate \eqref{S3 Linfty convergence} is weaker 
compared to the results in the Euclidean setting, where one has
\begin{align*}
    t^{\frac{n}{2}}\,
    \|u(t,\,\cdot\,)-Mh_{t}\|_{L^{\infty}(\mathbb{R}^{n})}\,
    \longrightarrow\,0
    \qquad\textnormal{as}\quad\,t\rightarrow\infty.
\end{align*}
Such a convergence is no longer valid in our setting due to the spectral gap.
A counterexample was provided on $\mathbb{H}^{3}(\mathbb{R})$ in
\cite[p.15]{Vaz2019}.
Here, let us extend it to general noncompact symmetric spaces.
Consider a delayed heat kernel $h_{t+t'}$ for some $t'>0$ to be determined 
later. Then 
\begin{align*}
    t^{\frac{\nu}{2}}\,e^{|\rho|^{2}t}\,
    \|h_{t+t'}-h_{t}\|_{L^{\infty}(\mathbb{X})}\,
    \ge\,t^{\frac{\nu}{2}}\,e^{|\rho|^{2}t}\,\big(h_{t}(eK)-h_{t+t'}(eK)\big)
\end{align*}
since $h_{t}(eK)$ is decreasing in $t$.
According to \eqref{S2 heat kernel}, there exists a constant $C\ge1$ such that
\begin{align*}
    t^{\frac{\nu}{2}}\,e^{|\rho|^{2}t}\,\big(h_{t}(eK)-h_{t+t'}(eK)\big)
    &\ge\,
    t^{\frac{\nu}{2}}\,e^{|\rho|^{2}t}\,
    \big\lbrace{
        \tfrac{1}{C}\,t^{-\frac{\nu}{2}}\,e^{-|\rho|^{2}t}\,
        -\,
        C(t+t')^{-\frac{\nu}{2}}\,e^{-|\rho|^{2}(t+t')}
    }\big\rbrace\\[5pt]
    &=\,
    C^{-1}\,-\,C\big(\tfrac{t}{t+t'}\big)^{\frac{\nu}{2}}\,e^{-|\rho|^{2}t'}
    \ge\,\tfrac{1}{2C}\,,
\end{align*}
provided that $t'>\frac{2\ln{C}+\ln2}{|\rho|^{2}}$.
Hence
\begin{align*}
    t^{\frac{\nu}{2}}\,e^{|\rho|^{2}t}\,
    \|h_{t+t'}-h_{t}\|_{L^{\infty}(\mathbb{X})}\,
     \centernot\longrightarrow\,0
    \qquad\textnormal{as}\quad\,t\rightarrow\infty.
\end{align*}
\end{remark}
\subsection{Counterexample in the non bi-$K$-invariant case}\label{S3 Sub4}
Recall that the heat kernel concentrates in the region $K(\exp\Omega_{t})K$
where $\Omega_{t}=B(2t\rho,r(t))$. Here $r(t)$ goes to infinity slightly 
faster than $\sqrt{t}$ and slower than $t$, as $t\rightarrow\infty$. 
Fix $y\in{G}\smallsetminus{K}$ and
consider the heat equation with the initial data $u_{0}(x)=\delta_{y}(x)$,
where $\delta_{y}$ denotes the Dirac measure at $y$. Since $\delta_{y}$ is not bi-K-invariant, it is sufficient to show that the property \eqref{S1 Main thm 1 convergence} breaks down with such initial data. We write
\begin{align}
    \|u(t,\,\cdot\,)-Mh_{t}\|_{L^{1}(\mathbb{X})}\,
    &=\,\|h_{t}(\,\cdot\,,yK)
        -h_{t}(\,\cdot\,,eK))\|_{L^{1}(\mathbb{X})}\notag\\[5pt]
    &\ge\,
    \Big\|h_{t}(\,\cdot\,,eK)
    \Big(\frac{h_{t}(\,\cdot\,,yK)}{h_{t}(\,\cdot\,,eK)}-1\Big)
    \Big\|_{L^{1}(K(\exp\Omega_{t})K)}\label{S3 counterexample difference of two kernels}\\[5pt]
    &=\,
    \Big\lbrace{
    \int_{\Omega_{t}}\diff{g^{+}}\,\delta(g^{+})\,h_{t}(\exp{g^{+}})\,
    }\Big\rbrace\,
    \Big\lbrace{
    \int_{K}\diff{k}\,|e^{\langle{2\rho,\,A(k^{-1}y)}\rangle}-1|
    +\textrm{O}\big(\tfrac{r(t)}{t}\big)}\Big\rbrace
\notag\end{align}
according to the Cartan decomposition 
$g=k(\exp{g^{+}})k'$ and \cref{S3 proposition kernel quotient} which we will prove in a moment.
Since the heat kernel $h_{t}$ is bi-$K$-invariant and concentrates in
$K(\exp\Omega_{t})K$, we have
\begin{align*}
    \int_{\Omega_{t}}\diff{g^{+}}\,\delta(g^{+})\,h_{t}(\exp{g^{+}})\,
    =\,\|h_{t}\|_{L^{1}(K(\exp\Omega_{t})K)}\,
    \longrightarrow\,1
    \qquad\textnormal{as}\quad\,t\rightarrow\infty.
\end{align*}
Hence, the right-hand side of 
\eqref{S3 counterexample difference of two kernels} tends to
\begin{align*}
    \int_{K}\diff{k}\,|e^{\langle{2\rho,\,A(k^{-1}y)}\rangle}-1|
\end{align*}
as $t\rightarrow\infty$, which is not identically $0$. 
Indeed, if $y=k_{1}(\exp{y^{+}})k_{2}$ in the Cartan decomposition,
then the nonnegative continuous function
$k\mapsto|e^{\langle{2\rho,\,A(k^{-1}y)}\rangle}-1|$
does not vanish at $k=k_{1}$, where 
$\langle{\rho,\,A(k^{-1}y)}\rangle=\langle{\rho,\,y^{+}}\rangle>0$,
since $y\notin{K}$ implies $y^{+}\neq0$.
In conclusion, \cref{S1 Main thm 1} fails without the 
bi-$K$-invariance assumption. 
It remains to prove \cref{S3 proposition kernel quotient}.

\begin{proposition}\label{S3 proposition kernel quotient}
Let $y\in{G}\smallsetminus{K}$.
Then, for every $g$ in the critical region $K(\exp\Omega_{t})K$, we have
\begin{align*}
    \frac{h_{t}(gK,yK)}{h_{t}(gK,eK)}\,
    =\,e^{\langle{2\rho,\,A(k^{-1}y)}\rangle}
        +\textnormal{O}\big(\tfrac{r(t)}{t}\big)
    \qquad\textnormal{as}\quad\,t\rightarrow\infty,
\end{align*}
where $g=k(\exp{g^{+}})k'$ in the Cartan decomposition.
\end{proposition}
For the proof of Proposition \ref{S3 proposition kernel quotient} we need the following lemma.

\begin{lemma}\label{S3 lemma distance behaviors}
For all $g$ in the critical region $K(\exp\Omega_{t})K$,
the following asymptotic behaviors hold as $t\rightarrow\infty$: 
\begin{enumerate}[label=\textnormal{(\roman*)}]
    \vspace{5pt}\item
        $\frac{|(y^{-1}g)^{+}|}{|g^{+}|}
            =1+\textnormal{O}\big(\tfrac{r(t)}{t}\big)$.
    
    \vspace{5pt}\item
        $\frac{g^{+}}{|g^{+}|}$ and $\frac{(y^{-1}g)^{+}}{|(y^{-1}g)^{+}|}$
        are both equal to
        $\frac{\rho}{|\rho|}+\textnormal{O}\big(\tfrac{r(t)}{t}\big)$. 
        
    \vspace{5pt}\item
        For every $\alpha\in\Sigma^{+}$,
        $\frac{\langle{\alpha,(y^{-1}g)^{+}}\rangle}{
            \langle{\alpha,g^{+}}\rangle}
            =1+\textnormal{O}\big(\tfrac{r(t)}{t}\big)$.
    
    \vspace{5pt}\item
        $d(gK,eK)-d(gK,yK)
            =\langle{\frac{\rho}{|\rho|},A(k^{-1}y)}\rangle
                +\textnormal{O}\big(\tfrac{r(t)}{t}\big)$.
\end{enumerate}
\end{lemma}

\begin{proof}[Proof of \cref{S3 lemma distance behaviors}]
For all $g\in{K(\exp{B(2t\rho,r(t))})K}$ and for any fixed $y\in{G}$,
we have
\begin{align*}
    d(gK,eK)\,=\,2|\rho|t+\textrm{O}\big(r(t)\big)
    \qquad\textnormal{and}\qquad
    d(gK,yK)\,=\,2|\rho|t+\textrm{O}\big(r(t)\big)
\end{align*}
as $t\rightarrow\infty$. 
We deduce first (i) and (ii) by using
\begin{align*}
    \frac{|(y^{-1}g)^{+}|}{|g^{+}|}\,
    =\,\frac{d(gK,yK)}{d(gK,eK)}\,
    =\,1+\textrm{O}\big(\tfrac{r(t)}{t}\big),
    \qquad
    \frac{g^{+}}{|g^{+}|}
    =\,
    \frac{2t\rho+\textrm{O}\big(r(t)\big)}{2t|\rho|+\textrm{O}\big(r(t)\big)}\,
    =\,\frac{\rho}{|\rho|}+\textrm{O}\big(\tfrac{r(t)}{t}\big),    
\end{align*}
and 
\begin{align*}
    |(y^{-1}g)^{+}-g^{+}|\,
    =\,|(g^{-1}y)^{+}-(g^{-1})^{+}|
    \le\,d(g^{-1}yK,g^{-1}K)\,
    =\,|y|\,=\,\mathrm{O}(1).
\end{align*}
Here we have used \eqref{S2 Distance} and the fact that
\begin{align*}
    (x^{-1})^{+}\,=\,-w_{0}.\,x^{+}
    \qquad\forall\,x\in{G},
\end{align*}
where $w_{0}$ denotes the longest element in the Weyl group $W$,
which exchanges the positive and the negative Weyl chambers.
Let us next deduce (iii) from (i) and (ii). 
For every positive root $\alpha$,
\begin{align*}
    \frac{\langle{\alpha,(y^{-1}g)^{+}}\rangle}{\langle{\alpha,g^{+}}\rangle}\,
    &=\,
    \frac{\langle{\alpha,\frac{(y^{-1}g)^{+}}{|(y^{-1}g)^{+}|}}\rangle}{
        \langle{\alpha,\frac{g^{+}}{|g^{+}|}}\rangle}\,
        \frac{|(y^{-1}g)^{+}|}{|g^{+}|}\notag\\[5pt]
    &=\, 
    \frac{\langle{\alpha,\frac{\rho}{|\rho|}}\rangle
        +\textrm{O}\big(\tfrac{r(t)}{t}\big)}{
        \langle{\alpha,\frac{\rho}{|\rho|}}\rangle
        +\textrm{O}\big(\tfrac{r(t)}{t}\big)}\,
        \Big\lbrace{1+\textrm{O}\big(\tfrac{r(t)}{t}\big)}\Big\rbrace\,
    =\,1+\textrm{O}\big(\tfrac{r(t)}{t}\big).
\end{align*}
It remains us to prove (iv). Let $g=k(\exp{g^{+}})k'$
in the Cartan decomposition and consider the Iwasawa decomposition
$k^{-1}y=n(k^{-1}y)(\exp{A(k^{-1}y)})k''$ for some $k''\in{K}$. Then
\begin{align}\label{S3 distance decomposition in lemma}
    d(gK,yK)\,
    &=\,d\big(k(\exp{g^{+}})K,kn(k^{-1}y)(\exp{A(k^{-1}y)})K\big)
        \notag\\[5pt]
    &=\,d\big(\exp{(-g^{+})}[n(k^{-1}y)]^{-1}(\exp{g^{+}})K,
        \exp{(A(k^{-1}y)}-g^{+})K\big).
\end{align}
and we write
\begin{align*}
    d(gK,eK)-d(gK,yK)\,
    &=\,\overbrace{\vphantom{\Big|}
        d(gK,eK)-d\big(\exp{(A(k^{-1}y)}-g^{+})K,eK\big)}^{I}\\
    &+\,\underbrace{\vphantom{\Big|}
        d\big(\exp{(A(k^{-1}y)}-g^{+})K,eK\big)-d(gK,yK)}_{II}.
\end{align*}
On the one hand, we have
\begin{align*}
    I\,
    =\,|g^{+}|-|A(k^{-1}y)-g^{+}|\,
    &=\,\frac{2\langle{g^{+},A(k^{-1}y)}\rangle-|A(k^{-1}y)|^{2}}{
        |g^{+}|+|A(k^{-1}y)-g^{+}|}\\[5pt]
    &=\,\big\langle{\tfrac{g^{+}}{|g^{+}|},A(k^{-1}y)}\big\rangle\,
        +\textrm{O}\big(\tfrac{1}{|g^{+}|}\big)\\[5pt]
    &=\,\big\langle{\tfrac{\rho}{|\rho|},A(k^{-1}y)}\big\rangle\,
        +\textrm{O}\big(\tfrac{r(t)}{t}\big)
\end{align*}
by using $(ii)$ and the fact that $\lbrace{A(k^{-1}y)\,|\,k\in{K}}\rbrace$ 
is a compact subset of $\mathfrak{a}$.
On the other hand, we deduce from \eqref{S3 distance decomposition in lemma} 
that
\begin{align*}
    |II|\,
    \le\,d\big(\exp{(-g^{+})}[n(k^{-1}y)]^{-1}(\exp{g^{+}})K,eK\big)
\end{align*}
where the right-hand side tends exponentially fast to $0$, as
$\lbrace{n(k^{-1}y)\,|\,k\in{K}}\rbrace$ is a compact 
subset of $N$ see for instance \cite[Lemma 3.1]{Hel1963}.
Indeed, the subgroup $\exp\mathfrak{a}$ acts by conjugation on
$\mathfrak{n}$ and $N$. In particular $\exp(-g^+)$ acts by contractions
$e^{-\langle\alpha,g^+\rangle}$ on each subspace $\mathfrak{g}_\alpha$.
Thus $\Ad(\exp(-g^+))(\log{n})$ tends to $0$ exponentially fast, for every
$\log{n}\in\mathfrak{n}$ and more generally uniformly on compact subsets of
$\mathfrak{n}$. Via the exponential map, $\exp(-g^+)n(\exp g^+)$ tends to $e$
exponentially fast, for every $n\in{N}$ and more generally uniformly 
on compact subsets of $N$.
In conclusion, we obtain
\begin{align}\label{S3 Busemann function}
    d(gK,eK)-d(gK,yK)\,
    =\,\langle{\tfrac{\rho}{|\rho|},A(k^{-1}y)}\rangle
        +\textrm{O}\big(\tfrac{r(t)}{t}\big).
\end{align}
\end{proof}

\begin{remark}
Observe that \eqref{S3 Busemann function} defines in the limit
of a Busemann function. Let us elaborate.
Any regular geodesic ray in $G/K$ stemming from the origin is given by
$\gamma(r)=k(\exp{rH_{0}})K$, where $k\in{K}$ and $H_{0}$ is a unit vector in 
$\mathfrak{a}^{+}$. 
As in the proof of $(iv)$ in \cref{S3 lemma distance behaviors},
consider the Iwasawa decomposition $k^{-1}y=n(k^{-1}y)(\exp{A(k^{-1}y)})k''$
for a given $y\in{G}$. We have, on the one hand $d(\gamma(r),eK)=r$.
On the other hand,
\begin{align*}
    d(\gamma(r),yK)\,
    =\,d\big(\exp{(-rH_{0})}[n(k^{-1}y)]^{-1}(\exp{rH_{0}})K,
        \exp{(A(k^{-1}y)}-rH_{0})K\big).
\end{align*}
Hence
\begin{align*}
    |d(\gamma(r),yK)-
    \underbrace{\vphantom{\Big|}
        d(\exp{(A(k^{-1}y)-rH_{0}})K,eK)|}_{|A(k^{-1}y)-rH_{0}|}\,
    \le\,d\big(
        \underbrace{\vphantom{\Big|}
        \exp{(-rH_{0})}[n(k^{-1}y)]^{-1}(\exp{rH_{0}})}_{\rightarrow\,e}
        K,eK\big)
\end{align*}
as $r\rightarrow\infty$. Therefore,
\begin{align*}
    d(\gamma(r),yK)-d(\gamma(r),eK)\,
    &=\,|A(k^{-1}y)-rH_{0}|-r+\textnormal{o}(1)\\[5pt]
    &=\,\frac{|A(k^{-1}y)|^{2}-2r\langle{H_{0},A(k^{-1}y)}\rangle}{
        2r+\textnormal{O}(1)}+\textnormal{o}(1)\\[5pt]
    &=\,-\langle{H_{0},A(k^{-1}y)}\rangle+\textnormal{o}(1)
\end{align*}
as $r\rightarrow\infty$. In conclusion, 
we have thus determined the Busemann function
\begin{align*}
    B_{\gamma}(yK)\,
    =\,\lim_{r\rightarrow+\infty}
        \lbrace{d(\gamma(r),yK)-d(\gamma(r),eK)}\rbrace
    =\,-\langle{H_{0},A(k^{-1}y)}\rangle.
\end{align*}
\end{remark}

Now, let us turn to the proof of \cref{S3 proposition kernel quotient}.
\begin{proof}[Proof of \cref{S3 proposition kernel quotient}]
Using \eqref{S2 heat kernel critical region} and
\eqref{S2 phi0 far}, we write
\begin{align*}
    \frac{h_{t}(gK,yK)}{h_{t}(gK,eK)}\,
    &=\,\frac{h_{t}(\exp(y^{-1}g)^{+})}{h_{t}(\exp(g^{+}))}\\[5pt]
    &\sim\,
        \frac{\mathbf{b}\big(-i\tfrac{(y^{-1}g)^{+}}{2t}\big)^{-1}}{
            \mathbf{b}\big(-i\tfrac{g^{+}}{2t}\big)^{-1}}\,
        \frac{\bm{\pi}((y^{-1}g)^{+})}{\bm{\pi}(g^{+})}\,
        \frac{e^{-\langle{\rho,(y^{-1}g)^{+}}\rangle}}{
            e^{-\langle{\rho,g^{+}}\rangle}}\,
        \frac{e^{-\frac{|(y^{-1}g)^{+}|^{2}}{4t}}}{
            e^{-\frac{|g^{+}|^{2}}{4t}}}
\end{align*}
as $t\rightarrow\infty$. 
The asymptotic behaviors of the first two factors are based on
\cref{S3 lemma distance behaviors}.(iii).
On the one hand,
\begin{align}\label{S3 counterexample 2}
    \frac{\bm{\pi}((y^{-1}g)^{+})}{\bm{\pi}(g^{+})}\,
    =\,\prod_{\alpha\in\Sigma_{r}^{+}}
    \frac{\langle{\alpha,(y^{-1}g)^{+}}\rangle}{\langle{\alpha,g^{+}}\rangle}\,
    =\,1+\textrm{O}\big(\tfrac{r(t)}{t}\big).
\end{align}
On the other hand, by using \eqref{S2 behavior of Gamma functions}, we have
\begin{align}\label{S3 counterexample 1}
    \frac{\mathbf{b}\big(-i\tfrac{(y^{-1}g)^{+}}{2t}\big)^{-1}}{
            \mathbf{b}\big(-i\tfrac{g^{+}}{2t}\big)^{-1}}\,
    \sim\,\prod_{\alpha\in\Sigma_{r}^{+}}
    \Big\lbrace
    \frac{\langle{\alpha,g^{+}}\rangle}{\langle{\alpha,(y^{-1}g)^{+}}\rangle}
    \Big\rbrace^{1-\frac{m_{\alpha}}{2}-\frac{m_{2\alpha}}{2}}\,
    =\,1+\textrm{O}\big(\tfrac{r(t)}{t}\big).
\end{align}
For the last two factors, we notice on the one hand that
\begin{align}\label{S3 norm distance}
    -\frac{|(y^{-1}g)^{+}|^{2}}{4t}+\frac{|g^{+}|^{2}}{4t}\,
    &=\,\frac{|g^{+}|+|(y^{-1}g)^{+}|}{4t}\,
        \big(|g^{+}|-|(y^{-1}g)^{+}|\big)\notag\\[5pt]
    &=\,\underbrace{\vphantom{\Big|}
        \frac{d(gK,eK)+d(gK,yK)}{4t}
        }_{=\,|\rho|+\textrm{O}\big(\tfrac{r(t)}{t}\big)}
        \underbrace{\vphantom{\frac{|g^{+}|+|(y^{-1}g)^{+}|}{4t}}
        \big\lbrace{d(gK,eK)-d(gK,yK)}\big\rbrace
        }_{=\,\langle{\frac{\rho}{|\rho|},A(k^{-1}y)}\rangle
                +\textrm{O}\big(\tfrac{r(t)}{t}\big)}\notag\\[5pt]
    &=\,\langle{\rho,A(k^{-1}y)}\rangle
                +\textrm{O}\big(\tfrac{r(t)}{t}\big)
\end{align}
according to \cref{S3 lemma distance behaviors}. Therefore
\begin{align}\label{S3 counterexample 3}
    e^{-\frac{|(y^{-1}g)^{+}|^{2}}{4t}+\frac{|g^{+}|^{2}}{4t}}\,
    =\,e^{\langle{\rho,A(k^{-1}y)}\rangle}+\textrm{O}\big(\tfrac{r(t)}{t}\big).
\end{align}
On the other hand, since $(g^{-1})^{+}=-w_{0}.g^{+}$ and $\rho=-w_{0}.\rho$,
where $w_{0}\in{W}$ is the longest element in the Weyl group, we write
\begin{align*}
    -\langle{\rho,(y^{-1}g)^{+}}\rangle+\langle{\rho,g^{+}}\rangle\,
    &=\,
    -\langle{\rho,(g^{-1}y)^{+}}\rangle+\langle{\rho,(g^{-1})^{+}}\rangle\\[5pt]
    &=\,
    \underbrace{\vphantom{\Big|}
    \frac{|(g^{-1})^{+}|^{2}-|(g^{-1}y)^{+}|^{2}}{4t}
    }_{=\,\langle{\rho,A(k^{-1}y)}\rangle+\textrm{O}\big(\tfrac{r(t)}{t}\big)}
    +
    \underbrace{\vphantom{\Big|}
    \frac{|(g^{-1}y)^{+}-2t\rho|^{2}-|(g^{-1})^{+}-2t\rho|^{2}}{4t}}_{R}.
\end{align*}
where the first part of the right-hand side is estimated as 
\eqref{S3 norm distance}, while the remainder
\begin{align*}
    R\,=\,
    \big\lbrace
    |(g^{-1}y)^{+}-2t\rho|-|(g^{-1})^{+}-2t\rho|
    \big\rbrace\,
    \frac{|(g^{-1}y)^{+}-2t\rho|+|(g^{-1})^{+}-2t\rho|}{4t}
\end{align*}
is $\mathrm{O}(\tfrac{r(t)}{t}\big)$,
according to (ii) in \cref{S3 lemma distance behaviors}. Hence
\begin{align}\label{S3 counterexample 4}
    e^{-\langle{\rho,(y^{-1}g)^{+}}\rangle+\langle{\rho,g^{+}}\rangle}\,
    =\,e^{\langle{\rho,A(k^{-1}y)}\rangle}+\textrm{O}\big(\tfrac{r(t)}{t}\big).
\end{align}
In conclusion, we deduce from \eqref{S3 counterexample 1}, 
\eqref{S3 counterexample 2}, \eqref{S3 counterexample 3} and 
\eqref{S3 counterexample 4} that 
\begin{align*}
    \frac{h_{t}(gK,yK)}{h_{t}(gK,eK)}\,
    =\,e^{\langle{2\rho,\,A(k^{-1}y)}\rangle}
        +\textrm{O}\big(\tfrac{r(t)}{t}\big)
\end{align*}
with $\frac{r(t)}{t}\rightarrow0$ as $t\rightarrow\infty$.
\end{proof}

\begin{remark}\label{Remark 3.10} Let us notice that
\cref{S3 proposition kernel quotient} also  clarifies 
the role of bi-$K$-invariance, and in addition gives an alternative proof for \cref{S3 L1 convergence in critical region}.
Let $u_{0}\in C_c(\mathbb{X})$ and recall that $\mathbb{M}$ is the centralizer of $\exp{\mathfrak{a}}$ in $K$. 
Recall that the Helgason-Fourier transform 
    \begin{align}\label{S3 Helgason}
        \mathcal{H}u_{0}(\lambda,k\mathbb{M})\,
        =\,
        \int_{G}\diff{g}\,
        u_{0}(gK)\,e^{\langle{-i\lambda+\rho,\,A(k^{-1}g)}\rangle}
    \end{align} 
boils down to the transform \eqref{S2 HC transform}
when $u_{0}$ is bi-$K$-invariant. 
According to \cref{S3 proposition kernel quotient}, we know that, for all $g$ in the critical region $K(\exp\Omega_{t})K$,
\begin{align}
    u(t,gK)\,-\,Mh_{t}(gK)\,
    &=\,
    \int_{G}\diff{y}\,u_{0}(y)\,
    (h_{t}(y^{-1}g)-h_{t}(gK))\notag\\[5pt]
    &=\,
    h_{t}(gK)\,\int_{G}\diff{y}\,u_{0}(y)\,
    \Big(
    e^{\langle{2\rho,\,A(k^{-1}y)}\rangle}\,-\,1\,
    +\textrm{O}\big(\tfrac{r(t)}{t}\big)
    \Big)\notag\\[5pt]
    &=\,
    h_{t}(gK)\,
    \big(
    \mathcal{H}u_{0}(i\rho,k\mathbb{M})\,
    -\,\mathcal{H}u_{0}(-i\rho,k\mathbb{M})\,
    +\,\textrm{O}\big(\tfrac{r(t)}{t}\big)\,
    \mathcal{H}u_{0}(-i\rho,k\mathbb{M})
    \big).
    \label{S3 new rmk}
\end{align}
Notice that $\mathcal{H}u_{0}(\pm\,i\rho,k\mathbb{M})=\mathcal{H}u_{0}(\pm\,i\rho)=M$ when $u_{0}$ is bi-$K$-invariant. Then we deduce \cref{S3 L1 convergence in critical region} by integrating \eqref{S3 new rmk} over the critical region. On the other hand, we have
\begin{align*}
    \int_{K(\exp\Omega_{t})K}\diff{g}\,
    |u(t,gK)\,-\,Mh_{t}(gK)|\,
    &\longrightarrow\,
    \int_{K}\diff{k}\,
    \Big|
    \int_{G}\diff{y}\,u_{0}(y)\,
    \big(
    e^{\langle{2\rho,\,A(k^{-1}y)}\rangle}\,-\,1
    \big)
    \Big|\\[5pt]
    &=\,
    \int_{K}\diff{k}\,
    |\mathcal{H}u_{0}(i\rho,k\mathbb{M})\,
    -\,\mathcal{H}u_{0}(-i\rho,k\mathbb{M})|,
\end{align*}
as $t\rightarrow\infty$. The last integral is not constantly zero when $u_{0}$ is not bi-$K$-invariant. For example, as we have seen in \eqref{S3 counterexample difference of two kernels}, if $u_{0}$ is a Dirac measure supported on some point outside of $K$, the last integral does not vanish.
\end{remark}

\section{Asymptotic convergence associated with 
the distinguished Laplacian}\label{Section.4 Distinguished}

Let $S=N(\exp{\mathfrak{a}})=(\exp{\mathfrak{a}})N$ be the solvable group
occurring in the Iwasawa decomposition $G=N(\exp{\mathfrak{a}})K$. 
Then $S$ is identifiable, as a manifold, with the symmetric space 
$\mathbb{X}=G/K$. The distinguished Laplacian $\widetilde{\Delta}$ on $S$ is 
given by the conjugation of the shifted Laplace-Beltrami operator
$\Delta+|\rho|^{2}$ on $\mathbb{X}$:
\begin{align}\label{S4 dist laplacian}
    \widetilde{\Delta}\,
    =\,\widetilde{\delta}^{\frac12}\circ(\Delta+|\rho|^{2})
        \circ\widetilde{\delta}^{-\frac12}
\end{align}
where the modular function $\widetilde{\delta}$ of $S$ is defined by
\begin{align*}
    \widetilde{\delta}(g)\,
    =\,\widetilde{\delta}(n(\exp{A}))\,
    =\,e^{-2\langle{\rho,A}\rangle}
    \qquad\forall\,g\in{S}.
\end{align*}
Here $n=n(g)$ and $A=A(g)$ denotes respectively the $N$-component and the
$\mathfrak{a}$-component of $g$ in the Iwasawa decomposition.
The distinguished Laplacian $\widetilde{\Delta}$ is left-$S$-invariant and
self-adjoint with respect to the right-invariant Haar measure on $S$:
\begin{align*}
    \int_{S}\textrm{d}_{r}{g}\,f(g)
    =\,\int_{N}\diff{n}\int_{\mathfrak{a}}\diff{A}\,f(n(\exp{A}))\,
    =\,\int_{\mathfrak{a}}\diff{A}\,e^{2\langle{\rho,A}\rangle}
        \int_{N}\diff{n}f((\exp{A})n).
\end{align*}
Bear in mind the following different relations between the measures on $S$
and the unimodular Haar measure on $G$:
\begin{align}\label{S4 drdldg}
    \int_{S}\textrm{d}_{r}{g}\,f(g)\,
    =\,\int_{G}\diff{g}\,e^{2\langle{\rho,A(g)}\rangle}f(g)
    \qquad\textnormal{and}\qquad
    \int_{S}\textrm{d}_{\ell}{g}\,f(g)\,
    =\,\int_{G}\diff{g}\,f(g).
\end{align}

This section aims to study the asymptotic behavior of solutions to
the following Cauchy problem associated with the distinguished Laplacian:
\begin{align}\label{HE S}
    \partial_{t}\widetilde{v}(t,g)\,
    =\,\widetilde{\Delta}_{g}\widetilde{v}(t,g),
    \qquad
    \widetilde{v}(0,g)\,=\,\widetilde{v}_{0}(g),
\end{align}
where the corresponding heat kernel is given by 
$\widetilde{h}_{t}=\widetilde{\delta}^{\frac12}\,e^{|\rho|^{2}t}h_{t}$
in the sense that
\begin{align*}
    (e^{t\widetilde{\Delta}\,}f)(g)\,
    =\,(f*\widetilde{h}_{t})(g)\,
    =\,\int_{S}\textrm{d}_{\ell}{y}\,f(y)\,\widetilde{h}_{t}(y^{-1}g)\,
    =\,\int_{S}\textrm{d}_{r}{y}\,f(gy^{-1})\,\widetilde{h}_{t}(y).
\end{align*}
Here, we still denote by $*$ the convolution product on $S$ or on $G$. We refer 
to \cite{Bou1983,CGGM1991} for more details about the distinguished Laplacian.

\begin{remark}
Notice that $\widetilde{h}_{t}(g)\textrm{d}_{r}{g}$ is a probability measure
on $S$. Indeed, we know that the Abel transform of $e^{|\rho|^{2}t}h_{t}$ is 
the heat kernel on $\mathfrak{a}$, i.e.,
\begin{align*}
    \mathcal{A}\big(e^{|\rho|^{2}t}h_{t})(A)\,
    =\,e^{-\langle{\rho,A}\rangle}
        \int_{N}\diff{n}\,e^{|\rho|^{2}t}h_{t}(n(\exp{A}))\,
    =\,(4\pi{t})^{-\frac{\ell}{2}}e^{-\frac{|A|^{2}}{4t}},
\end{align*}
hence
\begin{align*}
    \int_{N}\diff{n}\,\widetilde{h}_{t}(n(\exp{A}))\,
    =\,(4\pi{t})^{-\frac{\ell}{2}}e^{-\frac{|A|^{2}}{4t}}
\end{align*}
and
\begin{align*}
    \int_{S}\textrm{d}_{r}{g}\,\widetilde{h}_{t}(g)\,
    =\,\int_{N}\diff{n}\int_{\mathfrak{a}}\diff{A}\,
        \widetilde{h}_{t}(n(\exp{A}))\,
    =\,(4\pi{t})^{-\frac{\ell}{2}}
        \int_{\mathfrak{a}}\diff{A}\,e^{-\frac{|A|^{2}}{4t}}\,
    =\,1.
\end{align*}
\end{remark}

The first subsection is devoted to determine the critical region 
where the heat kernel $\widetilde{h}_{t}$ concentrates. 
In the next two subsections, we study respectively the $L^{1}$ 
and the $L^{\infty}$ asymptotic convergences of solutions to \eqref{HE S}
with compactly supported initial data (no bi-$K$-invariance required). 
We discuss the same questions for other initial data in the last subsection.

\subsection{Asymptotic concentration of the distinguished heat kernel}
Recall that the heat kernel $h_{t}$
(associated with the heat equation \eqref{HE X}) concentrates in 
$K(\exp\Omega_{t})K$ where $\Omega_{t}$ is a ball in the positive Weyl chamber.
The following proposition shows that the heat kernel $\widetilde{h}_{t}$ 
(associated with the Cauchy problem \eqref{HE S}) concentrates in a different region.
Recall that the function
$\mu:\overline{\mathfrak{a}^{+}}\rightarrow\mathbb{R}^{+}$
is defined by $\mu(H)=\min_{\alpha\in\Sigma^{+}}\langle{\alpha,H}\rangle$.

\begin{proposition}
Let $t\mapsto\varepsilon(t)$ be a positive function such that 
$\varepsilon(t)\searrow0$ and $\varepsilon(t)\sqrt{t}\rightarrow\infty$ 
as $t\rightarrow\infty$. 
Then the heat kernel associated with the distinguished Laplacian on $S$
concentrates asymptotically in $K(\exp\widetilde{\Omega}_{t})K$, where
\begin{align*}
    \widetilde{\Omega}_{t}\,
    =\,\big\lbrace{
        H\in\overline{\mathfrak{a}^{+}}\,|\,
        \varepsilon(t)\sqrt{t}\le|H|\le\tfrac{\sqrt{t}}{\varepsilon(t)}
        \,\,\,\textnormal{and}\,\,\,
        \mu(H)\ge\varepsilon(t)\sqrt{t}
    }\big\rbrace.
\end{align*}
In other words,
\begin{align*}
    \lim_{t\rightarrow\infty}
    \int_{g\in{S}\,\textrm{s.t.}\,
        g^{+}\in\overline{\mathfrak{a}^{+}}
        \smallsetminus\widetilde{\Omega}_{t}}
    \textrm{d}_{r}{g}\,\widetilde{h}_{t}(g)\,
    =\,0
\end{align*}
where $g^{+}$ denotes the middle component of $g$ in the Cartan decomposition.
\end{proposition}

\begin{figure}
    \centering
    \input{Concentration2.tex}
    \caption{Critical region $\widetilde{\Omega}_{t}$ 
        in the positive Weyl chamber.}
    \label{fig concentration2}
\end{figure}
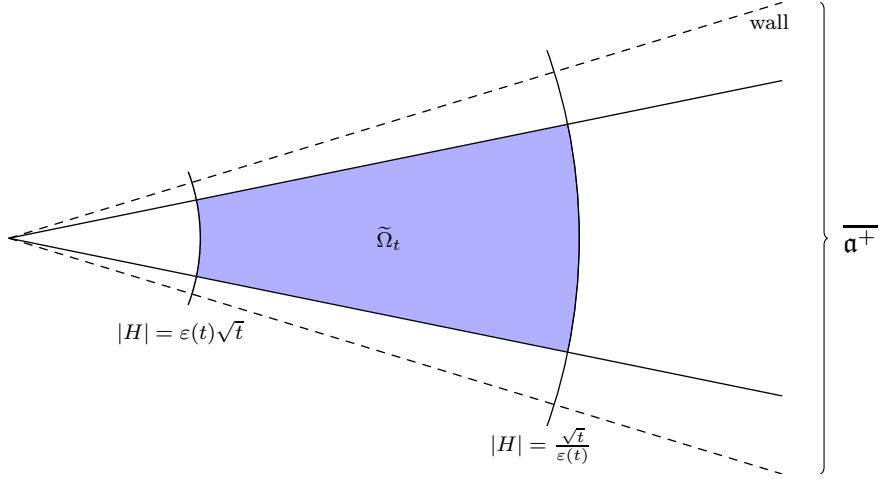

\begin{proof}
By using \eqref{S4 drdldg}, let us write
\begin{align*}
    I(t)\,
    =\,\int_{S\cap{K(\exp\widetilde{\Omega}_{t})K}}
        \textrm{d}_{r}{g}\,\widetilde{h}_{t}(g)\,
    =\,e^{|\rho|^{2}t}\int_{K(\exp\widetilde{\Omega}_{t})K}\diff{g}\,
        e^{\langle{\rho,A(g)}\rangle}\,h_{t}(g).
\end{align*}
Since the heat kernel $h_{t}$ is bi-$K$-invariant on $G$ and 
$\diff{k}$ is the normalized Haar measure on the compact group $K$, we have
\begin{align*}
    I(t)\,
    =\,e^{|\rho|^{2}t}
        \int_{K(\exp\widetilde{\Omega}_{t})K}\diff{g}\,h_{t}(g)
        \int_{K}\diff{k}\,e^{\langle{\rho,A(kg)}\rangle}\,
    =\,e^{|\rho|^{2}t}
        \int_{K(\exp\widetilde{\Omega}_{t})K}\diff{g}\,h_{t}(g)\,
        \varphi_{0}(g).
\end{align*}
According to the Cartan decomposition, and to the estimates
\eqref{S2 estimate of delta}, \eqref{S2 global estimate phi0} 
and \eqref{S2 heat kernel}, we obtain
\begin{align}\label{S4 global estimate I}
    I(t)\,
    &=\,e^{|\rho|^{2}t}
        \int_{\widetilde{\Omega}_{t}}\diff{g^{+}}\,
        \delta(g^{+})\,h_{t}(\exp{g^{+}})\,\varphi_{0}(\exp{g^{+}})
        \notag\\[5pt]
    &\asymp\,t^{-\frac{n}{2}}
        \int_{\widetilde{\Omega}_{t}}\diff{g^{+}}\,
        e^{-\frac{|g^{+}|^{2}}{4t}}\,
        \omega_{1}(g^{+})\,\omega_{2}(t,g^{+}) 
\end{align}
where
\begin{align*}
    \omega_{1}(g^{+})\,
    =\,\prod_{\alpha\in\Sigma^{+}}
        \Big( 
        \frac{\langle\alpha,g^{+}\rangle}
        {1+\langle\alpha,g^{+}\rangle}
        \Big)^{m_{\alpha}}\,
    \le\,1
\end{align*}
and
\begin{align*}
    \omega_{2}(t,g^{+})\,
    =\,\prod_{\alpha\in\Sigma_{r}^{+}}
        (1+\langle{\alpha,g^{+}}\rangle)^{2}
        (1+t+\langle{\alpha,g^{+}}\rangle)^{
        \frac{m_{\alpha}+m_{2\alpha}}{2}-1}.
\end{align*}
Next, let us study the right-hand side of \eqref{S4 global estimate I} 
outside $\widetilde{\Omega}_{t}$. On the one hand, let
\begin{align*}
    I_{0}(t)\,
    =\,t^{-\frac{n}{2}}
        \int_{|g^{+}|<\varepsilon(t)\sqrt{t}}\diff{g^{+}}\,
        e^{-\frac{|g^{+}|^{2}}{4t}}\,\omega_{2}(t,g^{+}).
\end{align*}
By substituting $g^{+}=2\sqrt{t}H$ and noticing that
\begin{align*}
    \omega_{2}(t,2\sqrt{t}H)\,
    &=\,(4t)^{\frac{n-\ell}{2}}
        \prod_{\alpha\in\Sigma_{r}^{+}}
        \big(\tfrac{1}{2\sqrt{t}}+\langle{\alpha,H}\rangle
        \big)^{2}
        \big(\tfrac{1}{2t}+\tfrac{1}{2}
        +\tfrac{\langle{\alpha,H}\rangle}{\sqrt{t}}
        \big)^{\frac{m_{\alpha}+m_{2\alpha}}{2}-1}\\[5pt]
    &\lesssim\,t^{\frac{n-\ell}{2}}
        \prod_{\alpha\in\Sigma_{r}^{+}}
        \big(1+\langle{\alpha,H}\rangle
        \big)^{\frac{m_{\alpha}+m_{2\alpha}+3}{2}}
\end{align*}
we get the upper bound of $I_{0}(t)$ in large time:
\begin{align}\label{S4 upper bound I0}
    I_{0}(t)\,
    \lesssim\,\int_{|H|<\frac{\varepsilon(t)}{2}}\diff{H}\,
        e^{-|H|^{2}}\,
        \prod_{\alpha\in\Sigma_{r}^{+}}
        \big(1+\langle{\alpha,H}\rangle
        \big)^{\frac{m_{\alpha}+m_{2\alpha}+3}{2}}
    \lesssim\,\varepsilon(t)^{\ell}.
\end{align}
On the other hand, similar computations yield
\begin{align}\label{S4 upper bound Iinf}
    I_{\infty}(t)\,
    &=\,t^{-\frac{n}{2}}
        \int_{|g^{+}|>\frac{\sqrt{t}}{\varepsilon(t)}}\diff{g^{+}}\,
        e^{-\frac{|g^{+}|^{2}}{4t}}\,\omega_{2}(t,g^{+})\notag\\[5pt]
    &\lesssim\,
        \int_{|H|>\frac{1}{2\varepsilon(t)}}\diff{H}\,
        e^{-|H|^{2}}\,
        |H|^{\tfrac{n-\ell+3|\Sigma_{r}^{+}|}{2}}\,
    \lesssim\,\varepsilon(t)^{N}
\end{align}
for any $N\ge0$. In the case where 
$\mu(g^{+})=\min_{\alpha\in\Sigma^{+}}\langle{\alpha,g^{+}}\rangle
<\varepsilon(t)\sqrt{t}$, i.e., when $g^{+}$ is close to the walls,
we have 
\begin{align*}
    I_{\mu}(t)\,
    =\,t^{-\frac{n}{2}}
        \int_{\mu(g^{+})<\varepsilon(t)\sqrt{t}}\diff{g^{+}}\,
        e^{-\frac{|g^{+}|^{2}}{4t}}\,\omega_{2}(t,g^{+})\,
    =\,2^{\ell}\,t^{-\frac{n-\ell}{2}}
        \int_{\mu(H)<\frac{\varepsilon(t)}{2}}\diff{H}\,
        e^{-|H|^{2}}\,\omega_{2}(t,2\sqrt{t}H)
\end{align*}
after substituting $g^{+}=2\sqrt{t}H$. 
Notice that, there exists at least one $\alpha_{0}\in\Sigma_{r}^{+}$ 
such that $\langle{\alpha_{0},H}\rangle<\frac{\varepsilon(t)}{2}$.
For such an $\alpha_{0}$, we estimate
\begin{align*}
    \tfrac{1}{2\sqrt{t}}+\langle{\alpha_{0},H}\rangle\,
    <\,\tfrac{1}{2\sqrt{t}}+\tfrac{\varepsilon(t)}{2}\,
    \lesssim\,\varepsilon(t)\,(1+\langle{\alpha_{0},H}\rangle).
\end{align*}
For the other $\alpha\in\Sigma_{r}^{+}\smallsetminus\lbrace\alpha_{0}\rbrace$, 
we simply estimate 
\begin{align*}
    \tfrac{1}{2\sqrt{t}}+\langle{\alpha,H}\rangle
    \lesssim\,\big(1+\langle{\alpha,H}\rangle\big)^{2}.
\end{align*}
Altogether, 
\begin{align}\label{S4 upper bound Im}
    I_{\mu}(t)\,
    \lesssim\,
        \varepsilon(t)^{2}
        \underbrace{\vphantom{\Big|}
        \int_{\mu(H)<\frac{\varepsilon(t)}{2}}\diff{H}\,
        e^{-|H|^{2}}\,
        \prod_{\alpha\in\Sigma_{r}^{+}}
        \big(1+\langle{\alpha,H}\rangle
        \big)^{\frac{m_{\alpha}+m_{2\alpha}+3}{2}}
        }_{<\,\infty}.
\end{align}
In conclusion, we deduce from \eqref{S4 global estimate I},
\eqref{S4 upper bound I0}, \eqref{S4 upper bound Iinf} and 
\eqref{S4 upper bound Im}, that
\begin{align*}
    0\,\le\,
    \int_{g\in{S}\,\textrm{s.t.}\,
        g^{+}\in\overline{\mathfrak{a}^{+}}
        \smallsetminus\widetilde{\Omega}_{t}}
    \textrm{d}_{r}{g}\,\widetilde{h}_{t}(g)\,
    \lesssim\,\varepsilon(t)
\end{align*}
for $t$ large enough. Since $\varepsilon(t)\rightarrow0$ as $t\rightarrow\infty$,
we have proved that the heat kernel $\widetilde{h}_{t}$ associated with the 
distinguished Laplacian on $S$ concentrates asymptotically in 
$K(\exp\widetilde{\Omega}_{t})K$.
\end{proof}

The asymptotic behaviors \eqref{S2 heat kernel critical region} and
\eqref{S2 phi0 far} of the heat kernel $h_{t}(\exp{H})$ and the ground 
spherical function $\varphi_{0}(\exp{H})$ hold when 
$\mu(H)\rightarrow\infty$, i.e., when $H\in\overline{\mathfrak{a}^{+}}$ stays 
away from the walls. The following proposition sharpens these asymptotics
when $H$ lies in the present critical region $\widetilde{\Omega}_{t}$. 

\begin{proposition}
Let $C_{1},C_{2}$ be the positive constants occurring in 
\eqref{S2 heat kernel critical region}, \eqref{S2 phi0 far},
and let $C_{3}=C_{1}\bm{\pi}(\rho_{0})^{-1}$.
Then the following asymptotics hold for 
$H\in\widetilde{\Omega}_{t}$ and $t\rightarrow\infty$:
\begin{align}\label{S4 asymp ht}
    h_{t}(\exp{H})\,
    =\,\Big\lbrace{C_{3}+\textnormal{O}
        \big(\tfrac{1}{\varepsilon(t)\sqrt{t}}}\big)
        \Big\rbrace\,
        t^{-\frac{\nu}{2}}\,e^{-|\rho|^{2}t}\,
        \bm{\pi}(H)\,e^{-\langle{\rho,H}\rangle}\,
        e^{-\frac{|H|^2}{4t}}
\end{align}
and
\begin{align}\label{S4 asymp phi0}
    \varphi_{0}(\exp{H})\,
    =\,\Big\lbrace{C_{2}+\textnormal{O}
        \big(\tfrac{1}{\mu(H)}}\big)\Big\rbrace\,
        \bm{\pi}(H)\,e^{-\langle{\rho,H}\rangle}.
\end{align}
\end{proposition}

\begin{proof}
As $H\in\widetilde{\Omega}_{t}$ stays away from the walls in 
$\overline{\mathfrak{a}^{+}}$, we substitute the spherical function 
occurring in
\begin{align*}
    e^{|\rho|^{2}t}\,h_{t}(\exp{H})\,
    =\,C_{0}\int_{\mathfrak{a}}\diff{\lambda}\,
        |\mathbf{c}(\lambda)|^{-2}\,e^{-t|\lambda|^{2}}\,
        \varphi_{\lambda}(\exp{H})
\end{align*}
by its Harish-Chandra expansion \eqref{S2 Harish-Chandra expansion 1},
and we obtain firstly
\begin{align*}
    e^{|\rho|^{2}t}\,e^{\langle{\rho,H}\rangle}\,h_{t}(\exp{H})\,
    =\,C_{0}|W|\,
        \sum_{q\in2Q}e^{-\langle{q,H}\rangle}\,
        \int_{\mathfrak{a}}\diff{\lambda}\,
        \mathbf{c}(-\lambda)^{-1}\,e^{-t|\lambda|^{2}}\,
        e^{i\langle{\lambda,H}\rangle}\,\gamma_{q}(\lambda).
\end{align*}
By factorizing 
$\mathbf{c}(-\lambda)^{-1}=\mathbf{b}(-\lambda)^{-1}\bm{\pi}(-i\lambda)$
and by performing an integration by parts based on 
\begin{align*}
    \bm{\pi}(-i\lambda)e^{-t|\lambda|^{2}}\,
    =\,\bm{\pi}(\tfrac{i}{2t}\tfrac{\partial}{\partial\lambda})
        e^{-t|\lambda|^{2}}
\end{align*}
we obtain secondly 
\begin{align*}
    &(2t)^{|\Sigma_{r}^{+}|}\,e^{|\rho|^{2}t}\,
    e^{\langle{\rho,H}\rangle}\,h_{t}(\exp{H})\\[5pt]
    &=\,C_{0}|W|\,
        \sum_{q\in2Q}e^{-\langle{q,H}\rangle}\,
        \int_{\mathfrak{a}}\diff{\lambda}\,e^{-t|\lambda|^{2}}\,
        \bm{\pi}(-i\tfrac{\partial}{\partial\lambda})\,
        \lbrace{
        e^{i\langle{\lambda,H}\rangle}
        \mathbf{b}(-\lambda)^{-1}\gamma_{q}(\lambda)
        }\rbrace\\[5pt]
    &=\,C_{0}|W|\,
        \sum_{q\in2Q}e^{-\langle{q,H}\rangle}\,
        \sum_{\Sigma_{r}^{+}=\Sigma_{1}\sqcup\Sigma_{2}\sqcup\Sigma_{3}}
        \int_{\mathfrak{a}}\diff{\lambda}\,e^{-t|\lambda|^{2}}\\[5pt]
        &\times\,\Big\lbrace{e^{i\langle{\lambda,H}\rangle}
        \prod_{\alpha_{1}\in\Sigma_{1}}\langle{\alpha_{1},H}\rangle
        }\Big\rbrace \,
        \Big\lbrace{
        \prod_{\alpha_{2}\in\Sigma_{2}}(-i\partial_{\alpha_2})
        (\mathbf{b}(-\lambda)^{-1})
        }\Big\rbrace
        \Big\lbrace{
        \prod_{\alpha_{3}\in\Sigma_{3}}(-i\partial_{\alpha_3})
        (\gamma_{q}(\lambda))
        }\Big\rbrace.
\end{align*}
After the shift of contour $\lambda\mapsto\lambda+i\frac{H}{2t}$ 
in $\mathfrak{a}_{\mathbb{C}}$ and the rescaling 
$\lambda\mapsto\frac{\lambda}{\sqrt{t}}$, we obtain thirdly
\begin{align}\label{S4 third integral}
    &t^{\frac{\nu}{2}}\,e^{|\rho|^{2}t}\,e^{\langle{\rho,H}\rangle}\,
    e^{\frac{|H|^2}{4t}}h_{t}(\exp{H})\notag\\[5pt]
    &=\,C_{0}2^{-|\Sigma_{r}^{+}|}|W|\,
        \sum_{q\in2Q}e^{-\langle{q,H}\rangle}\,
        \sum_{\Sigma_{r}^{+}=\Sigma_{1}\sqcup\Sigma_{2}\sqcup\Sigma_{3}}
        \Big\lbrace{
        \prod_{\alpha_{1}\in\Sigma_{1}}\langle{\alpha_{1},H}\rangle
        }\Big\rbrace
        \int_{\mathfrak{a}}\diff{\lambda}\,e^{-|\lambda|^{2}}\notag\\[5pt]
        &\times\,\Big\lbrace{
        \prod_{\alpha_{2}\in\Sigma_{2}}
            (-i\partial_{\alpha_2})\,\mathbf{b}^{-1}
        }\Big\rbrace\,(-\tfrac{\lambda}{\sqrt{t}}-i\tfrac{H}{2t})\,
        \Big\lbrace{
        \prod_{\alpha_{3}\in
            \Sigma_{3}}(-i\partial_{\alpha_3})\,\gamma_{q}
        }\Big\rbrace
        \,(-\tfrac{\lambda}{\sqrt{t}}-i\tfrac{H}{2t}).
\end{align}
Recall that, for all $H\in\widetilde{\Omega}_{t}$, we have
$\varepsilon(t)\sqrt{t}\le|H|\le\frac{\sqrt{t}}{\varepsilon(t)}$
and $\mu(H)=\min_{\alpha\in\Sigma^{+}}\langle{\alpha,H}\rangle
\ge\varepsilon(t)\sqrt{t}$.
Under these assumptions, we notice that the derivatives of
$\mathbf{b}^{-1}$ are bounded by
\begin{align}\label{S4 estimate b}
    \big|\mathbf{b}(-\tfrac{\lambda}{\sqrt{t}}-i\tfrac{H}{2t})^{-1}\big|\,
    &\asymp\,
    \prod_{\alpha\in\Sigma_{r}^{+}}
        \big(1+\tfrac{\langle{\alpha,\lambda}\rangle}{\sqrt{t}}
        +\tfrac{\langle{\alpha,H}\rangle}{t}
        \big)^{\frac{m_{\alpha}+m_{2\alpha}}{2}-1}\notag\\[5pt]
    &\lesssim\,\prod_{\alpha\in\Sigma_{r}^{+}}
        (1+\langle{\alpha,\lambda}\rangle)^{\frac{m_{\alpha}+m_{2\alpha}}{2}}\,
    \lesssim\,(1+|\lambda|)^{\frac{n-\ell}{2}},
\end{align}
according to \eqref{S2 bfunction} and \eqref{S2 bfunction derivative},
while the derivatives of $\gamma_{q}$ are bounded by
\begin{align}\label{S4 estimate gammaq}
    |\gamma_{q}(-\tfrac{\lambda}{\sqrt{t}}-i\tfrac{H}{2t})|\,
    \lesssim\,(1+|q|)^{N_{\gamma}},
\end{align}
for some nonnegative constant $N_{\gamma}$, 
according to \eqref{S2 gammaq 1} and \eqref{S2 gammaq 2}. 
Let us split up the the right-hand side of \eqref{S4 third integral} as
\begin{align*}
    I(t,H)\,+\,R(t,H)
\end{align*}
where the leading term
\begin{align*}
    I(t,H)\,
    =\,C_{0}2^{-|\Sigma_{r}^{+}|}|W|\,\bm{\pi}(H)\,
        \int_{\mathfrak{a}}\diff{\lambda}\,e^{-|\lambda|^{2}}\,
        \mathbf{b}(-\tfrac{\lambda}{\sqrt{t}}-i\tfrac{H}{2t})^{-1}
\end{align*}
is the contribution of $q=0$ and  $\Sigma_{1}=\Sigma_{r}^{+}$,
while $R(t,H)$ denotes the remainder. 
On the one hand, noticing from \eqref{S4 estimate b} that 
\begin{align*}
    \big|\mathbf{b}(-\tfrac{\lambda}{\sqrt{t}}-i\tfrac{H}{2t})^{-1}
    -\mathbf{b}(0)^{-1}\big|\,
    \lesssim\,(1+|\lambda|)^{\frac{n-\ell}{2}}
        \big(\tfrac{|\lambda|}{\sqrt{t}}+\tfrac{|H|}{2t}\big)\,
    \lesssim\,\tfrac{1}{\varepsilon(t)\sqrt{t}}\,
        (1+|\lambda|)^{\frac{n-\ell}{2}+1},
\end{align*}
we obtain
\begin{align}\label{S4 asymp case2}
    I(t,H)\,=\,
    \big\lbrace{
    C_{3}+\textrm{O}\big(\tfrac{1}{\varepsilon(t)\sqrt{t}}\big)
    }\big\rbrace\,\bm{\pi}(H)
\end{align}
where $C_{3}=C_{0}2^{-|\Sigma_{r}^{+}|}|W|\pi^{\frac{\ell}{2}}\mathbf{b}(0)^{-1}
=C_{1}\bm{\pi}(\rho_{0})^{-1}$. On the other hand, we deduce from 
\eqref{S4 estimate b} and \eqref{S4 estimate gammaq} that
\begin{align*}
    |R(t,H)|\,
    \lesssim\,\bm{\pi}(H)\,
        \underbrace{\vphantom{\Bigg|}
        \sum\nolimits_{q\in2Q,\,q\neq0}e^{-\langle{q,H}\rangle}\,
            (1+|q|)^{N_{\gamma}}
        }_{\lesssim\,e^{-\mu(H)}}\,
        \underbrace{\vphantom{\Bigg|}
        \int_{\mathfrak{a}}\diff{\lambda}\,e^{-|\lambda|^{2}}\,
        (1+|\lambda|)^{\frac{n-\ell}{2}}
        }_{<\,\infty}.
\end{align*}
This proves the asymptotic behavior \eqref{S4 asymp ht}.
Let us next turn to \eqref{S4 asymp phi0}. Along the lines of \cite{AnJi1999},
we multiply \eqref{S2 Harish-Chandra expansion 1} by $\bm{\pi}(i\lambda)$ and
obtain
\begin{align}\label{S4 first equality phi0}
    e^{\langle{\rho,H}\rangle}\,\bm{\pi}(i\lambda)\,
    \varphi_{\lambda}(\exp{H})\,
    &=\,\sum_{w\in{W}}(\det{w})\,\mathbf{b}(w.\lambda)\,
        \Phi_{w.\lambda}(H)\notag\\[5pt]
    &=\,\sum\nolimits_{q\in2Q}e^{-\langle{q,H}\rangle}\,
        \sum_{w\in{W}}(\det{w})\,\mathbf{b}(w.\lambda)\,
        \gamma_{q}(w.\lambda)\,e^{i\langle{w.\lambda,H}\rangle}.
\end{align}
Here we used the factorization 
$\mathbf{c}(w.\lambda)=\bm{\pi}(iw.\lambda)^{-1}\mathbf{b}(w.\lambda)
=(\det{w})\bm{\pi}(i\lambda)^{-1}\mathbf{b}(w.\lambda)$ for every $w\in{W}$.
After applying $\bm{\pi}(-i\frac{\partial}{\partial\lambda})|_{\lambda=0}$ to
\eqref{S4 first equality phi0}, the left-hand side becomes
\begin{align*}
    |W|\,\bm{\pi}(\rho_{0})\,
    e^{\langle{\rho,H}\rangle}\,\varphi_{0}(\exp{H})
\end{align*}
(see for instance the beginning of the proof of Proposition 2.2.12 in
\cite{AnJi1999}).
The leading term 
\begin{align*}
    I'(t,H)\,
    =\,|W|\,\mathbf{b}(0)\,\bm{\pi}(H)
\end{align*}
in the right-hand side of \eqref{S4 first equality phi0} is obtained
by taking $q=0$ and by applying 
$\bm{\pi}(-i\frac{\partial}{\partial\lambda})|_{\lambda=0}$
to $e^{i\langle{w.\lambda,H}\rangle}$, while the remainder $R'(t,H)$ is 
estimated by 
\begin{align*}
    |R'(t,H)|\,
    \lesssim\,
    \big\lbrace{
    \tfrac{1}{\mu(H)}+e^{-\mu(H)}
    }\big\rbrace\,\bm{\pi}(H)\,
    \lesssim\,
    \tfrac{1}{\mu(H)}\,\bm{\pi}(H).
\end{align*}
This proves the asymptotic behavior \eqref{S4 asymp phi0}.
\end{proof}

\subsection{Heat asymptotics in $L^1$ for compactly supported initial data}
In this subsection, we investigate the long-time asymptotic convergence in 
$L^1(S)$ of solutions to the Cauchy problem \eqref{HE S}, where the initial data 
$\widetilde{v}_{0}$ is assumed continuous and compactly supported in $B(eK,\xi)$. 
For every $\lambda\in\mathfrak{a}$, let
$\widetilde{\varphi}_{\lambda}=\widetilde{\delta}^{\frac12}\varphi_{\lambda}$
be the modified spherical function. The mass function is defined by
\begin{align}\label{S4 mass G}
    \widetilde{M}(g)\,
    =\,\frac{(\widetilde{v}_{0}*\widetilde{\varphi}_{0})(g)
        }{\widetilde{\varphi}_{0}(g)}
    \qquad\forall\,g\in{S}.
\end{align}
By using the fact that the modular function $\widetilde{\delta}$ is a character
on $S$, we can also write the mass as
\begin{align}\label{S4 mass S}
    \widetilde{M}(g)\,
    =\,\tfrac{1}{\widetilde{\delta}(g)^{\frac12}\,\varphi_{0}(g)}\,
        \int_{S}\textrm{d}_{\ell}{y}\,
        v_{0}(gK)\,
        \underbrace{\vphantom{\Big|}
        \widetilde{\delta}(y)^{\frac12}
        \widetilde{\delta}(y^{-1}g)^{\frac12}
        }_{\widetilde{\delta}(g)^{\frac12}}\,
        \varphi_{0}(y^{-1}g)
    =\,\frac{(v_{0}*\varphi_{0})(g)}{\varphi_{0}(g)}
\end{align}
where $v_{0}(gK)=\widetilde{\delta}(g)^{-\frac{1}{2}}\widetilde{v}_{0}(g)$
is a right $K$-invariant function on $G$, 
with compact support $(\supp\widetilde{v}_{0})K$.

\begin{remark}\label{S4 mass remark}
The map
\begin{align*}
    \widetilde{v}_{0}\longmapsto\widetilde{v}_{0}*\widetilde{\varphi}_{0}
\end{align*}
can be interpreted as the spectral projection at the bottom $0$ of the spectrum
of $\widetilde{\Delta}$. Thus $\widetilde{M}$ generalizes somehow the mass in the
Euclidean case. Let us elaborate. We recall the Helgason-Fourier transform \eqref{S3 Helgason} and its inverse formula
\begin{align}\label{S4 Helgason inverse}
    f(gK)\,
    =\,|W|^{-1}\,
        \int_{\mathfrak{a}}\,\frac{\diff{\lambda}}{|\mathbf{c}(\lambda)|^{2}}\,
        \int_{K}\diff{k}\,\mathcal{H}f(\lambda,kM)\,
            e^{\langle{-i\lambda-\rho,H(g^{-1}k)}\rangle}.
\end{align}
By combining \eqref{S3 Helgason} with \eqref{S4 Helgason inverse} 
and by using the formula
\begin{align}\label{S4 spherical split}
    \varphi_{\lambda}(y^{-1}g)\,
    =\,\int_{K}\diff{k}\,e^{\langle{i\lambda+\rho,A(kg})\rangle}\,
        e^{\langle{-i\lambda+\rho,A(ky})\rangle}
\end{align}
(see for instance \cite[Chap.III, Theorem 1.1]{Hel1994}), we obtain
\begin{align*}
    v_{0}(gK)\,
    =\,|W|^{-1}\,\int_{\mathfrak{a}}\diff{\lambda}\,
        |\mathbf{c}(\lambda)|^{2}\,(v_{0}*\varphi_{\lambda})(g)
\end{align*}
with
\begin{align*}
    \Delta(v_{0}*\varphi_{\lambda})\,
    =\,v_{0}*(\Delta\varphi_{\lambda})\,
    =\,-(|\rho|^{2}+|\lambda|^{2})(v_{0}*\varphi_{\lambda}).
\end{align*}
Hence
\begin{align*}
    \widetilde{v}_{0}(g)\,
    =\,|W|^{-1}\,\int_{\mathfrak{a}}\diff{\lambda}\,
        |\mathbf{c}(\lambda)|^{2}\,
        (\widetilde{v}_{0}*\widetilde{\varphi}_{\lambda})(g)
\end{align*}
with
\begin{align*}
    \widetilde{\Delta}(\widetilde{v}_{0}*\widetilde{\varphi}_{\lambda})\,
    =\,\widetilde{v}_{0}*(\widetilde{\Delta}\widetilde{\varphi}_{\lambda})\,
    =\,-|\lambda|^{2}(\widetilde{v}_{0}*\widetilde{\varphi}_{\lambda}).
\end{align*}
\end{remark}

\begin{remark}
If $\widetilde{v}_{0}\in\mathcal{C}_{c}(S)$,
then the mass function $\widetilde{M}$ is bounded.
This follows indeed from the local Harnack inequality:
\begin{align}\label{S4 Harnack}
    \varphi_{0}(y^{-1}g)\,
    =\,
    \int_{K}\diff{k}\,e^{\langle{\rho,A(kg})\rangle}\,
        e^{\langle{\rho,A(ky)}\rangle}\,
    \lesssim\,\int_{K}\diff{k}\,e^{\langle{\rho,A(kg})\rangle}\,
    =\,\varphi_{0}(g)
\end{align}
which holds for every $y\in\supp{v_{0}}=(\supp{\widetilde{v}_{0}})K$ and for every
$g\in{G}$. Here, we have used \eqref{S4 spherical split} and the fact that
$|A(ky)|\le|y|$ is bounded. See also \cite[Proposition 4.6.3]{GaVa1988}.
\end{remark}

\begin{remark}
If $\widetilde{v}_{0}=\delta^{1/2}v_{0}$ with $v_{0}$ bi-$K$-invariant,
we notice on the one hand that
\begin{align}\label{S4 radial equiv L1}
    \int_{S}\textrm{d}_{r}{g}\,\delta(g)^{\frac12}v_{0}(g)\,
    &=\,\int_{G}\diff{g}\,v_{0}(g)\,e^{\langle{\rho,A(g)}\rangle}\notag\\[5pt]
    &=\,\int_{G}\diff{g}\,v_{0}(g)\,
        \int_{K}\diff{k}\,e^{\langle{\rho,A(kg)}\rangle}\,
    =\,\int_{G}\diff{g}\,v_{0}(g)\,\varphi_{0}(g)
\end{align}
and on the other hand that
\begin{align}\label{S4 radial mass}
    \frac{(v_{0}*\varphi_{0})(g)}{\varphi_{0}(g)}\,
    &=\,\tfrac{1}{\varphi_{0}(g)}\,
        \int_{G}\diff{y}\,v_{0}(y)\varphi_{0}(y^{-1}g)\notag\\[5pt]
    &=\,\tfrac{1}{\varphi_{0}(g)}\,
        \int_{G}\diff{y}\,v_{0}(y)\,
        \underbrace{\vphantom{\Big|}
        \int_{K}\diff{k}\,\varphi_{0}(y^{-1}kg)
        }_{\varphi_{0}(g)\,\varphi_{0}(y)}\,
    =\,\int_{G}\diff{y}\,v_{0}(y)\,\varphi_{0}(y).
\end{align}
Hence the mass function $\widetilde{M}$ is a constant 
if $v_{0}$ is bi-$K$-invariant and 
$\widetilde{v}_{0}=\widetilde{\delta}^{\frac12}v_{0}$ belongs to $L^{1}(S)$:
\begin{align*}
    \widetilde{M}\,
    =\,\int_{G}\diff{y}\,v_{0}(y)\,\varphi_{0}(y)
    =\,\mathcal{H}v_{0}(0).
\end{align*}
\end{remark}

The following lemma plays a key role in the proof of \cref{S1 Main thm 2}.
\begin{lemma}\label{S4 Lemma ratios difference}
For  bounded $y\in{G}$ and for all $g$ in the critical region
$K(\exp\widetilde{\Omega}_{t})K$, the following asymptotic behavior holds:
\begin{align*}
    \frac{h_{t}(g^{-1}y)}{h_{t}(g^{-1})}
        -\frac{\varphi_{0}(g^{-1}y)}{\varphi_{0}(g^{-1})}\,
    =\,\mathrm{O}\Big(\frac{1}{\varepsilon(t)\sqrt{t}}\Big)
    \qquad\textnormal{as}\,\,\,t\longrightarrow\infty.
\end{align*}
\end{lemma}

\begin{proof}
Assume that $|y|\le\xi$ for some positive constant $\xi$.
Recall that for every $H\in\widetilde{\Omega}_{t}$, we have 
$\varepsilon(t)\sqrt{t}\le|H|\le\frac{\sqrt{t}}{\varepsilon(t)}$
and $\mu(H)\ge\varepsilon(t)\sqrt{t}$, where $\varepsilon(t)\rightarrow0$
and $\varepsilon(t)\sqrt{t}\rightarrow\infty$ as $t\rightarrow\infty$.
Since we are interested in the asymptotic behavior when $t$ goes to infinity,
we can assume that $t$ is large enough such that
\begin{align*}
    \varepsilon(t)\,<\sqrt{2}\,
    \qquad\textnormal{and}\qquad
    \varepsilon(t)\sqrt{t}\,>\,
        2\xi(1+\max\nolimits_{\alpha\in\Sigma}|\alpha|).
\end{align*}
Notice first that 
\begin{align*}
    g\in{K(\exp\widetilde{\Omega}_{t})K}
    \qquad\Longleftrightarrow\qquad
    g^{-1}\in{K(\exp\widetilde{\Omega}_{t})K}
\end{align*}
as
\begin{align*}
    H\in\widetilde{\Omega}_{t}
    \qquad\Longleftrightarrow\qquad
    -w_{0}.H\in\widetilde{\Omega}_{t}
\end{align*}
where $w_{0}$ denotes the longest element in $W$, which interchanges the 
positive and negative roots. Notice next that
\begin{align*}
    |(g^{-1}y)^{+}-(g^{-1})^{+}|\,
    \le\,d(g^{-1}yK,g^{-1}K)\,
    =\,|y|\,<\,\xi
\end{align*}
according to \eqref{S2 Distance}. Then we deduce the following estimates:
\begin{align*}
\begin{cases}
    |(g^{-1}y)^{+}|\,\le\,|(g^{-1})^{+}|+\xi\,
        <\,\tfrac{\sqrt{t}}{\varepsilon(t)}+\tfrac{1}{2}\varepsilon(t)\sqrt{t}\,
        <\,\tfrac{2\sqrt{t}}{\varepsilon(t)},\\[5pt]
    |(g^{-1}y)^{+}|\,\ge\,|(g^{-1})^{+}|-\xi\,
        >\,\varepsilon(t)\sqrt{t}-\tfrac{1}{2}\varepsilon(t)\sqrt{t}
        >\,\tfrac{\varepsilon(t)\sqrt{t}}{2},\\[5pt]
    \langle{\alpha,(g^{-1}y)^{+}}\rangle\,
        >\,\langle{\alpha,(g^{-1})^{+}}\rangle-|\alpha|\xi\,
        >\,\varepsilon(t)\sqrt{t}-\tfrac{1}{2}\varepsilon(t)\sqrt{t}
        >\,\tfrac{\varepsilon(t)\sqrt{t}}{2}
        \qquad\forall\alpha\in\Sigma^{+}.
\end{cases}
\end{align*}
In other words, we obtain
\begin{align*}
    g^{-1}y\,\in\,K(\exp{\widetilde{\Omega}_{t}'})K
    \qquad\forall\,g\in{K(\exp\widetilde{\Omega}_{t})K},\,\,
    \forall\,|y|<\xi,
\end{align*}
where
\begin{align*}
    \widetilde{\Omega}_{t}'\,
    =\,\big\lbrace{
        H\in\overline{\mathfrak{a}^{+}}\,|\,
        \varepsilon'(t)\sqrt{t}\le|H|\le\tfrac{\sqrt{t}}{\varepsilon'(t)}
        \,\,\,\textnormal{and}\,\,\,
        \mu(H)\ge\varepsilon'(t)\sqrt{t}
    }\big\rbrace.
\end{align*}
and $\varepsilon'(t)=\frac{1}{2}\varepsilon(t)$ is still a decreasing function
satisfying $\varepsilon'(t)\rightarrow0$ and 
$\varepsilon'(t)\sqrt{t}\rightarrow\infty$ as $t\rightarrow\infty$.
Thus the asymptotics \eqref{S4 asymp ht} and \eqref{S4 asymp phi0} yield
\begin{align*}
    \frac{h_{t}(g^{-1}y)}{h_{t}(g^{-1})}\,
    =\,\underbrace{
    \tfrac{C_{3}+\mathrm{O}\big(\tfrac{1}{\varepsilon'(t)\sqrt{t}}\big)}{
        C_{3}+\mathrm{O}\big(\tfrac{1}{\varepsilon(t)\sqrt{t}}\big)}
    }_{1+\mathrm{O}\big(\tfrac{1}{\varepsilon(t)\sqrt{t}}\big)}\,
    \underbrace{
    \vphantom{\tfrac{C_{3}+\mathrm{O}\big(\tfrac{1}{\varepsilon'(t)\sqrt{t}}\big)
        }{C_{3}+\mathrm{O}\big(\tfrac{1}{\varepsilon(t)\sqrt{t}}\big)}}
    \frac{\bm{\pi}((g^{-1}y)^{+})}{\bm{\pi}((g^{-1})^{+})}\,
    e^{-\langle{\rho,(g^{-1}y)^{+}}\rangle+\langle{\rho,(g^{-1})^{+}}\rangle}
    }_{R(g,y)}\,
    e^{-\tfrac{|(g^{-1}y)^{+}|^{2}}{4t}+\tfrac{|(g^{-1})^{+}|^{2}}{4t}}
\end{align*}
and
\begin{align*}
    \frac{\varphi_{0}(g^{-1}y)}{\varphi_{0}(g^{-1})}\,
    =\,
    \tfrac{C_{2}+\mathrm{O}\big(\tfrac{1}{\varepsilon'(t)\sqrt{t}}\big)}{
        C_{2}+\mathrm{O}\big(\tfrac{1}{\varepsilon(t)\sqrt{t}}\big)}
    \,R(g,y)\,
    =\,\big\lbrace{
        1+\mathrm{O}\big(\tfrac{1}{\varepsilon(t)\sqrt{t}}\big)
    }\big\rbrace\,R(g,y).
\end{align*}
Hence
\begin{align*}
    \frac{h_{t}(g^{-1}y)}{h_{t}(g^{-1})}
        -\frac{\varphi_{0}(g^{-1}y)}{\varphi_{0}(g^{-1})}\,
    =\,\Big\lbrace{
        1+\mathrm{O}\big(\tfrac{1}{\varepsilon(t)\sqrt{t}}\big)
        }\Big\rbrace\,
        \Big\lbrace{
        e^{-\tfrac{|(g^{-1}y)^{+}|^{2}}{4t}+\tfrac{|(g^{-1})^{+}|^{2}}{4t}}-1
        }\Big\rbrace\,R(g,y)
\end{align*}
On the one hand, since 
$\frac{\bm{\pi}((g^{-1}y)^{+})}{\bm{\pi}((g^{-1})^{+})}$ and 
$e^{-\langle{\rho,(g^{-1}y)^{+}}\rangle+\langle{\rho,(g^{-1})^{+}}\rangle}$
are uniformly bounded when $g\in{K(\exp{\widetilde{\Omega}_{t}})K}$ and
$|y|<\xi$, so is $R(g,y)$. On the other hand,
\begin{align*}
    e^{-\tfrac{|(g^{-1}y)^{+}|^{2}}{4t}+\tfrac{|(g^{-1})^{+}|^{2}}{4t}}\,
    =\,e^{\tfrac{|(g^{-1})^{+}|-|(g^{-1}y)^{+}|}{4\sqrt{t}}
        \tfrac{|(g^{-1})^{+}|+|(g^{-1}y)^{+}|}{\sqrt{t}}}\,
    =\,e^{\mathrm{O}\big(\tfrac{1}{\varepsilon(t)\sqrt{t}}\big)}\,
    =\,1+\mathrm{O}\big(\tfrac{1}{\varepsilon(t)\sqrt{t}}\big).
\end{align*}
In conclusion,
\begin{align*}
    \frac{h_{t}(g^{-1}y)}{h_{t}(g^{-1})}
        -\frac{\varphi_{0}(g^{-1}y)}{\varphi_{0}(g^{-1})}\,
    =\,\mathrm{O}\Big(\frac{1}{\varepsilon(t)\sqrt{t}}\Big)
        \qquad\forall\,g\in{K(\exp\widetilde{\Omega}_{t})K},\,\,
        \forall\,|y|<\xi.
\end{align*}
\end{proof}

\vspace{5pt}
Now, let us prove the first part of \cref{S1 Main thm 2}.
\begin{proof}[Proof of \eqref{S1 L1 disting} in \cref{S1 Main thm 2}]
By using
\begin{align*}
    (\widetilde{v}_{0}*\widetilde{\varphi}_{\lambda})(g)\,
    =\,\int_{S}\textrm{d}_{\ell}{y}\,
        v_{0}(yK)
        \underbrace{\vphantom{\Big|}
        \widetilde{\delta}(y)^{\frac12}\,\widetilde{\delta}(y^{-1}g)^{\frac12}
        }_{\widetilde{\delta}(g)^{\frac12}}
        \varphi_{\lambda}(y^{-1}g)
    =\,\widetilde{\delta}(g)^{\frac{1}{2}}
        (v_{0}*\varphi_{\lambda})(gK),
\end{align*}
let us write the solution $\widetilde{v}$ to \eqref{HE S} as
\begin{align*}
    \widetilde{v}(t,g)\,
    =\,(\widetilde{v}_{0}*\widetilde{h}_{t})(g)\,
    =\,e^{|\rho|^{2}t}\,\widetilde{\delta}(g)^{\frac{1}{2}}\,
        (v_{0}*h_{t})(g)
\end{align*}
We aim to study the difference
\begin{align}\label{S4 difference}
    \widetilde{v}(t,g)-\widetilde{M}(g)\widetilde{h}_{t}(g)\,
    &=\,\widetilde{h}_{t}(g)\,
        \frac{(v_{0}*h_{t})(g)}{h_{t}(g)}
        \,-\,
        \widetilde{h}_{t}(g)\,
        \frac{(v_{0}*\varphi_{0})(g)}{\varphi_{0}(g)}
        \notag\\[5pt]
    &=\,\widetilde{h}_{t}(g)\,
        \int_{G}\diff{y}\,v_{0}(yK)\,
        \Big\lbrace{
        \frac{h_{t}(y^{-1}g)}{h_{t}(g)}
        -\frac{\varphi_{0}(y^{-1}g)}{\varphi_{0}(g)}
        }\Big\rbrace\notag\\[5pt]
    &=\,\widetilde{h}_{t}(g)\,
        \int_{G}\diff{y}\,v_{0}(yK)\,
        \Big\lbrace{
        \frac{h_{t}(g^{-1}y)}{h_{t}(g^{-1})}
        -\frac{\varphi_{0}(g^{-1}y)}{\varphi_{0}(g^{-1})}
        }\Big\rbrace
\end{align}
where the last expression is derived from the symmetries 
$h_{t}(x^{-1})=h_{t}(x)$ and $\varphi_{0}(x^{-1})=\varphi_{0}(x)$.
According to the previous lemma, we have
\begin{align*}
    \frac{h_{t}(g^{-1}y)}{h_{t}(g^{-1})}
        -\frac{\varphi_{0}(g^{-1}y)}{\varphi_{0}(g^{-1})}\,
    =\,\mathrm{O}\Big(\frac{1}{\varepsilon(t)\sqrt{t}}\Big)
    \qquad\forall\,g\in{K(\exp\widetilde{\Omega}_{t})K},\,\,
    \forall\,y\in\supp{v_{0},}
\end{align*}
and therefore the integral of
$\widetilde{v}(t,\cdot)-\widetilde{M}\,\widetilde{h}_{t}$ 
over the critical region
\begin{align*}
    \int_{S\cap{K(\exp\widetilde{\Omega}_{t})K}}
        \textrm{d}_{r}{g}\,    
        |\widetilde{v}(t,g)-\widetilde{M}(g)\widetilde{h}_{t}(g)|\,
    \lesssim\,
        \frac{1}{\varepsilon(t)\sqrt{t}}\,
        \underbrace{\vphantom{\Big|}
        \int_{S}\textrm{d}_{r}{g}\,\widetilde{h}_{t}(g)
        }_{1}\,
        \underbrace{\vphantom{\Big|}
        \int_{G}\diff{y}\,|v_{0}(yK)|
        }_{\const}
\end{align*}
tends asymptotically to $0$. It remains for us to check that the integral
\begin{align*}
    \int_{S\smallsetminus{K(\exp\widetilde{\Omega}_{t})K}}\textrm{d}_{r}{g}\,
        |\widetilde{v}(t,g)-\widetilde{M}(g)\widetilde{h}_{t}(g)|\,
    &\le\,
    \int_{S\smallsetminus{K(\exp\widetilde{\Omega}_{t})K}}\textrm{d}_{r}{g}\,
        |\widetilde{v}(t,g)|\\[5pt]
    &+\,
    \int_{S\smallsetminus{K(\exp\widetilde{\Omega}_{t})K}}\textrm{d}_{r}{g}\,    
        |\widetilde{M}(g)|\widetilde{h}_{t}(g)
\end{align*}
tends also to $0$. On the one hand, we know that $\widetilde{M}$ is bounded
and that the heat kernel $\widetilde{h}_{t}$ asymptotically concentrates in 
$K(\exp\widetilde{\Omega}_{t})K$, hence
\begin{align*}
    \int_{S\smallsetminus{K(\exp\widetilde{\Omega}_{t})K}}\textrm{d}_{r}{g}\,
        |\widetilde{M}(g)|\widetilde{h}_{t}(g)\,
        \longrightarrow\,0
\end{align*}
as $t\rightarrow\infty$. On the other hand, notice that for all
$y\in\supp{v_{0}}$ and for all $g\in{G}$ such that
$g^{+}\notin\widetilde{\Omega}_{t}$, we have 
\begin{align}\label{S4 Omega''}
    (y^{-1}g)^{+}\,\notin\,\widetilde{\Omega}_{t}''
    =\,
    \big\lbrace{
        H\in\overline{\mathfrak{a}^{+}}\,|\,
        \varepsilon''(t)\sqrt{t}\le|H|\le\tfrac{\sqrt{t}}{\varepsilon''(t)}
        \,\,\,\textnormal{and}\,\,\,
        \mu(H)\ge\varepsilon''(t)\sqrt{t}
    }\big\rbrace
\end{align}
where $\varepsilon''(t)=2\varepsilon(t)$ (see the proof of the previous lemma).
Hence
\begin{align*}
    \int_{S\smallsetminus{K(\exp\widetilde{\Omega}_{t})K}}\textrm{d}_{r}{g}\,
        |\widetilde{v}(t,g)|\,
    &\le\,\int_{G}\diff{y}\,|v_{0}(yK)|\,
        \int_{G\smallsetminus{K(\exp\widetilde{\Omega}_{t})K}}\diff{g}\,
        \delta(g)^{-\frac12}\,e^{|\rho|^{2}t}\,h_{t}(y^{-1}g)\\[5pt]
    &\lesssim\,
        \underbrace{\vphantom{
        \int_{S\smallsetminus{K(\exp\widetilde{\Omega}_{t}'')K}}}
        \int_{S}\textrm{d}_{r}{y}\,|\widetilde{v}_{0}(y)|
        }_{\|\widetilde{v}_{0}\|_{L^{1}(S)}}\,
        \underbrace{
        \int_{S\smallsetminus{K(\exp\widetilde{\Omega}_{t}'')K}}
            \textrm{d}_{r}{g}\,\widetilde{h}_{t}(g)
        }_{\longrightarrow\,0}.\
\end{align*}
This concludes the proof of the heat asymptotics in $L^{1}$ for the 
distinguished Laplacian $\widetilde{\Delta}$ on $S$ and for initial data
$\widetilde{v}_{0}\in\mathcal{C}_{c}(S)$.
\end{proof}

\subsection{Heat asymptotics in $L^{\infty}$
for compactly supported initial data}

Let us start with the following lemma, which allows us to compare the middle
components occurring in the Iwasawa decomposition and in the Cartan decomposition.

\begin{lemma}
For all $g\in{G}$, we have
\begin{align}\label{S4 ineq Iwasawa Cartan}
    \langle{\rho,A(g)}\rangle\,
    \le\,\langle{\rho,g^{+}}\rangle
\end{align}
where $A(g)$ denotes the $\mathfrak{a}$-component of $g$ in the Iwasawa 
decomposition and $g^{+}$ denotes its $\overline{\mathfrak{a}^{+}}$-component
in the Cartan decomposition.
\end{lemma}

\begin{proof}
According to Kostant's convexity theorem \cite[Theorem IV.10.5]{Hel2000},
for every $H\in\mathfrak{a}$, the Iwasawa projection of $K(\exp{H})K$
is equal to the convex hull of $W.H$. 
Therefore $A(K(\exp{H})K)\subset\co(W.H)$, which implies that 
$A(g)\in\co(W.g^{+})$ for all $g\in{G}$.
The inequality \eqref{S4 ineq Iwasawa Cartan} follows from the fact that
$\rho\in\mathfrak{a}^{+}$ while $g^{+}-\co(W.g^{+})$ is contained in the
cone generated by the positive roots \cite[Lemma IV.8.3]{Hel2000}.
\end{proof}

In the following two propositions we collect some elementary properties of the 
distinguished heat kernel. The first one clarifies the lower
and the upper bounds of $\widetilde{h}_{t}$, while the second one describes
its critical region for the $L^{\infty}$ norm.

\begin{proposition}\label{S4 htilde estimate proposition}
The heat kernel $\widetilde{h}_{t}$ associated with the distinguished Laplacian
satisfies
\begin{align}\label{S4 htilde estimate}
    \|\widetilde{h}_{t}\|_{L^{\infty}(S)}\,
    \asymp\,t^{-\frac{\ell+|\Sigma_{r}^{+}|}{2}}
\end{align}
for $t$ large enough.
\end{proposition}

\begin{proof}
Using the global estimates \eqref{S2 global estimate phi0} and 
\eqref{S2 heat kernel}, we have
\begin{align}\label{S4 htilde1}
    \widetilde{h}_{t}(g)\,
    \asymp\,
    t^{-\frac{\ell+|\Sigma_{r}^{+}|}{2}}\,
    e^{-\langle{\rho,A(g)}\rangle}e^{-\langle{\rho,g^{+}}\rangle}\,
    e^{-\frac{|g^{+}|^2}{4t}}\,
    \Big\lbrace{
    \prod\nolimits_{\alpha\in\Sigma_{r}^{+}}
    \tfrac{1+\langle{\alpha,g^{+}}\rangle}{\sqrt{t}}
    \big(1+ \tfrac{\langle{\alpha,g^{+}}\rangle}{t}\big)^{
        \tfrac{m_{\alpha}+m_{2\alpha}}{2}-1}
    }\Big\rbrace
\end{align}
We obtain first the lower bound in \eqref{S4 htilde estimate}
by evaluating the right hand side of \eqref{S4 htilde1} at
$g_{0}=\exp(-\sqrt{t}\rho)$ and by observing that
\begin{align*}
    A(g_{0})\,=\,-\sqrt{t}\rho
    \qquad\textnormal{and}\qquad
    g_{0}^{+}\,=\,\sqrt{t}\rho.
\end{align*}
For the upper bound, notice that
\begin{align}\label{S4 htilde2}
     e^{-\langle{\rho,A(g)}\rangle}e^{-\langle{\rho,g^{+}}\rangle}\,
     \le\,1
\end{align}
according to \eqref{S4 ineq Iwasawa Cartan}, and that
\begin{align*}
    \tfrac{1+\langle{\alpha,g^{+}}\rangle}{\sqrt{t}}
    \big(1+ \tfrac{\langle{\alpha,g^{+}}\rangle}{t}\big)^{
        \tfrac{m_{\alpha}+m_{2\alpha}}{2}-1}\,
    \le\,
    \big(1+\tfrac{\langle{\alpha,g^{+}}\rangle}{\sqrt{t}}\big)
    \big(1+ \tfrac{\langle{\alpha,g^{+}}\rangle}{t}\big)^{
        \tfrac{m_{\alpha}+m_{2\alpha}}{2}}\,
    \le\,
    \big(1+ \tfrac{\langle{\alpha,g^{+}}\rangle}{\sqrt{t}}\big)^{
        \tfrac{m_{\alpha}+m_{2\alpha}}{2}+1}
\end{align*}
for $t$ large enough. We deduce from \eqref{S4 htilde1} that
\begin{align}\label{S4 htilde1''}
    \widetilde{h}_{t}(g)\,
    \lesssim\,t^{-\frac{\ell+|\Sigma_{r}^{+}|}{2}}\,
    \underbrace{
    e^{-\frac{1}{4}\big|\tfrac{g^{+}}{\sqrt{t}}\big|^{2}}
    \prod\nolimits_{\alpha\in\Sigma_{r}^{+}}
    \big(1+\langle{\alpha,\tfrac{g^{+}}{\sqrt{t}}}\rangle\big)^{
        \tfrac{m_{\alpha}+m_{2\alpha}}{2}+1}
    }_{=\,\textrm{O}(1)}\,
    \lesssim\,t^{-\frac{\ell+|\Sigma_{r}^{+}|}{2}}.
\end{align}
\end{proof}
\begin{proposition}
The heat kernel $\widetilde{h}_{t}$ concentrates asymptotically
in the same critical region for the $L^{\infty}$ norm as for the $L^1$ norm.
In other words,
\begin{align*}
    t^{\frac{\ell+|\Sigma_{r}^{+}|}{2}}\,
    \|\widetilde{h}_{t}\|_{L^{\infty}
    (S\smallsetminus{K(\exp\widetilde{\Omega}_{t})K})}\,
    \longrightarrow\,0
    \qquad\textnormal{as}\,\,\,t\rightarrow\infty.
\end{align*}
\end{proposition}

\begin{proof}
Let us study the sup norm of $\widetilde{h}_{t}$ outside the critical region. 
Recall that
\begin{align*}
    \widetilde{\Omega}_{t}\,
    =\,\big\lbrace{
        H\in\overline{\mathfrak{a}^{+}}\,|\,
        \varepsilon(t)\sqrt{t}\le|H|\le\tfrac{\sqrt{t}}{\varepsilon(t)}
        \,\,\,\textnormal{and}\,\,\,
        \mu(H)\ge\varepsilon(t)\sqrt{t}
    }\big\rbrace
\end{align*}
where $\mu(H)=\min_{\alpha\in\Sigma^{+}}\langle{\alpha,H}\rangle$, 
and $\varepsilon(t)\rightarrow0$ satisfies
$\tfrac{\sqrt{t}}{\varepsilon(t)}\rightarrow\infty$ and 
$\varepsilon(t)\sqrt{t}\rightarrow\infty$ as $t\rightarrow\infty$.
We deduce from  \eqref{S4 htilde1} and \eqref{S4 htilde2} that 
\begin{align}\label{S4 htilde1'}
    t^{\frac{\ell+|\Sigma_{r}^{+}|}{2}}\,
    \widetilde{h}_{t}(g)\,
    \lesssim\,
    e^{-\frac{|g^{+}|^2}{4t}}\,
    \prod\nolimits_{\alpha\in\Sigma_{r}^{+}}
    \tfrac{1+\langle{\alpha,g^{+}}\rangle}{\sqrt{t}}
    \big(1+ \tfrac{\langle{\alpha,g^{+}}\rangle}{t}\big)^{
        \tfrac{m_{\alpha}+m_{2\alpha}}{2}-1}
\end{align}
\noindent\textit{Case 1}: Assume that $|g^{+}|<\varepsilon(t)\sqrt{t}$.
Then we deduce easily from \eqref{S4 htilde1'} that
\begin{align*}
    t^{\frac{\ell+|\Sigma_{r}^{+}|}{2}}\,
    \widetilde{h}_{t}(g)\,
    \lesssim\,
    \varepsilon(t)^{|\Sigma_{r}^{+}|}
\end{align*}
tends to $0$.

\noindent\textit{Case 2}: Assume that $|g^{+}|>\tfrac{\sqrt{t}}{\varepsilon(t)}$.
The estimate \eqref{S4 htilde1''} implies that,
for any $N>0$,
\begin{align*}
    t^{\frac{\ell+|\Sigma_{r}^{+}|}{2}}\,
    \widetilde{h}_{t}(g)\,
    \lesssim\,
    e^{-\frac{1}{4}\big|\tfrac{g^{+}}{\sqrt{t}}\big|^{2}}
    \prod\nolimits_{\alpha\in\Sigma_{r}^{+}}
    \big(1+\langle{\alpha,\tfrac{g^{+}}{\sqrt{t}}}\rangle\big)^{
        \tfrac{m_{\alpha}+m_{2\alpha}}{2}+1}\,
    \lesssim\,\varepsilon(t)^{N},
\end{align*}
which tends to $0$.

\noindent\textit{Case 3}: Assume that $\mu(g^{+})<\varepsilon(t)\sqrt{t}$.
By using the estimate
\begin{align*}
    \tfrac{1+\langle{\alpha,g^{+}}\rangle}{\sqrt{t}}\,
    \lesssim\,\varepsilon(t)\,
    \lesssim\,\varepsilon(t)\,
        \big(\tfrac{1+\langle{\alpha,g^{+}}\rangle}{\sqrt{t}}\big)
\end{align*}
in \eqref{S4 htilde1'} for the roots $\alpha\in\Sigma_{r}^{+}$ satisfying
$\langle\alpha,g^{+}\rangle\le\varepsilon(t)\sqrt{t}$, we obtain
\begin{align*}
    t^{\frac{\ell+|\Sigma_{r}^{+}|}{2}}\,
    \widetilde{h}_{t}(g)\,
    \lesssim\,\varepsilon(t)\,
    \underbrace{
    e^{-\frac{1}{4}\big|\tfrac{g^{+}}{\sqrt{t}}\big|^{2}}
    \prod\nolimits_{\alpha\in\Sigma_{r}^{+}}
    \big(1+\langle{\alpha,\tfrac{g^{+}}{\sqrt{t}}}\rangle\big)^{
        \tfrac{m_{\alpha}+m_{2\alpha}}{2}+1}
    }_{=\,\textrm{O}(1)},
\end{align*}
which tends to $0$.

In conclusion, we obtain
\begin{align*}
    \widetilde{h}_{t}(g)\,=\,\mathrm{o}(t^{-\frac{\ell+|\Sigma_{r}^{+}|}{2}})
    \qquad\textnormal{as}\,\,\,t\rightarrow\infty
\end{align*}
outside the critical region $S\cap{K(\exp\widetilde{\Omega}_{t})K}$.
In other words, the heat kernel $\widetilde{h}_{t}(g)$ concentrates
asymptotically in the same critical region for the $L^{\infty}$ norm
as for the $L^{1}$ norm.
\end{proof}

\vspace{5pt}
Finally, let us prove the remaining part of \cref{S1 Main thm 2}.
\begin{proof}[Proof of \eqref{S1 Linf disting} in \cref{S1 Main thm 2}]
Let us resume the proof of \cref{S1 Main thm 2}.
In the critical region $S\cap{K(\exp\widetilde{\Omega}_{t})K}$, we have
\begin{align*}
    |\widetilde{v}(t,g)-\widetilde{M}(g)\widetilde{h}_{t}(g)|\,
    \le\,\widetilde{h}_{t}(g)\,
        \int_{|y|<\,\xi}\diff{g}\,|v_{0}(yK)|\,
        \Big|{
        \frac{h_{t}(g^{-1}y)}{h_{t}(g^{-1})}
        -\frac{\varphi_{0}(g^{-1}y)}{\varphi_{0}(g^{-1})}
        }\Big|
\end{align*}
with
\begin{align*}
    \Big|{
        \frac{h_{t}(g^{-1}y)}{h_{t}(g^{-1})}
        -\frac{\varphi_{0}(g^{-1}y)}{\varphi_{0}(g^{-1})}
    }\Big|\,
    \lesssim\,\frac{1}{\varepsilon(t)\sqrt{t}}
\end{align*}
according to \eqref{S4 difference} 
and to \cref{S4 Lemma ratios difference}.
Then we deduce from \eqref{S4 htilde estimate} that
\begin{align*}
    t^{\frac{\ell+|\Sigma_{r}^{+}|}{2}}\,
    |\widetilde{v}(t,g)-\widetilde{M}(g)\widetilde{h}_{t}(g)|\,
    \lesssim\,\tfrac{1}{\varepsilon(t)\sqrt{t}}
    \qquad\forall\,g\in{S\cap{K(\exp\widetilde{\Omega}_{t})K}}
\end{align*}
where the right-hand side tends to $0$ as $t\rightarrow\infty$.
Outside the critical region, we estimate separately $\widetilde{v}(t,g)$
and $\widetilde{M}(g)\widetilde{h}_{t}(g)$. On the one hand, we know that
$\widetilde{M}(g)$ is a bounded function and that 
$\widetilde{h}_{t}(g)=\mathrm{o}(t^{-(\ell+|\Sigma_{r}^{+}|)/2})$. Then 
$t^{(\ell+|\Sigma_{r}^{+}|)/2}\widetilde{M}(g)\widetilde{h}_{t}(g)$ 
tends to $0$ as $t\rightarrow\infty$.
On the other hand, since $g\notin{K(\exp\widetilde{\Omega}_{t})K}$ and 
$|y|<\xi$ imply that $g^{-1}y\notin{K(\exp\widetilde{\Omega}_{t}'')K}$
(see \eqref{S4 Omega''}), we obtain
\begin{align*}
    |\widetilde{v}(t,g)|\,
    \lesssim\,\int_{G}\diff{y}\,
        |\widetilde{v}_{0}(yK)|\,|\widetilde{h}_{t}(g^{-1}y)|\,
\end{align*}
which is $\textrm{o}(t^{-(\ell+|\Sigma_{r}^{+}|)/2})$ outside the critical region.
In conclusion,
\begin{align*}
    t^{\frac{\ell+|\Sigma_{r}^{+}|}{2}}
    \|\widetilde{v}(t,\,\cdot\,)-
    \widetilde{M}\,\widetilde{h}_{t}\|_{L^{\infty}(S)}\,
    \longrightarrow\,0
\end{align*}
as $t\rightarrow\infty$.
\end{proof}

By convexity we obtain easily the corresponding result for the $L^{p}$ norm.
\begin{corollary}
The solution $\widetilde{v}$ to the Cauchy problem \eqref{HE S} with initial 
data $\widetilde{v}_{0}\in\mathcal{C}_{c}(S)$ satisfies
\begin{align}\label{S4 Lp disting}
    t^{\frac{\ell+|\Sigma_{r}^{+}|}{2p'}}
    \|\widetilde{v}(t,\,\cdot\,)-
    \widetilde{M}\,\widetilde{h}_{t}\|_{L^{p}(S)}\,
    \longrightarrow\,0
    \qquad\textnormal{as}\quad\,t\rightarrow\infty,
\end{align}
for all $1<p<\infty$.
\end{corollary}

\subsection{Heat asymptotics for other initial data}\label{Subsect other data}
We have obtained above the long-time asymptotic convergence in $L^{p}$
($1\le{p}\le\infty$) for the distinguished heat equation 
with compactly supported initial data.
It is natural and interesting to ask whether the expected convergence
still holds when the initial data lie in larger functional spaces?
The following corollaries give positive examples,
but the optimal answer remains open.

\begin{corollary}
The asymptotic convergences \eqref{S1 L1 disting} and \eqref{S1 Linf disting}, 
hence \eqref{S4 Lp disting}, still hold with initial data
$\widetilde{v}_{0}=\widetilde{\delta}^{\tfrac12}v_{0}\in{L}^{1}(S)$ when 
$v_{0}$ is bi-$K$-invariant.
\end{corollary}

\begin{proof}
Notice from \eqref{S4 radial equiv L1} and \eqref{S4 radial mass} that 
the mass function $\widetilde{M}=(\mathcal{H}v_{0})(0)$ 
is a constant under the present assumption.
Let us start with the $L^{1}$ convergence. 
Given $\varepsilon>0$, let
$\widetilde{V}_{0}=\widetilde{\delta}^{\frac12}V_{0}\in\mathcal{C}_{c}(S)$,
with $V_{0}$ bi-$K$-invariant, be such that $\|\widetilde{v}_{0}
-\widetilde{V}_{0}\|_{L^{1}(S)}<\tfrac{\varepsilon}{3}$. 
Then the corresponding mass function $(\mathcal{H}V_{0})(0)$
is also a constant.
We observe firstly that the solution to the Cauchy problem
\begin{align*}
    \partial_{t}\widetilde{V}(t,g)\,
    =\,\widetilde{\Delta}_{g}\widetilde{V}(t,g),
    \qquad
    \widetilde{V}(0,g)\,=\,\widetilde{V}_{0}(g)
\end{align*}
satisfies
\begin{align}\label{S4 density1}
    \|\widetilde{v}(t,\,\cdot\,)-\widetilde{V}(t,\,\cdot\,)\|_{L^{1}(S)}\,
    =\,\|\widetilde{v}_{0}*\widetilde{h}_{t}
        -\widetilde{V}_{0}*\widetilde{h}_{t}\|_{L^{1}(S)}\,
    \le\,\|\widetilde{v}_{0}-\widetilde{V}_{0}\|_{L^{1}(S)}\,
        \|\widetilde{h}_{t}\|_{L^{1}(S)}\,
    <\,\tfrac{\varepsilon}{3},
\end{align}
since $\|\widetilde{h}_{t}\|_{L^{1}(S)}=1$, 
and secondly that there exists $T>0$ such that for all $t\ge{T}$,
\begin{align}\label{S4 density2}
    \|\widetilde{V}(t,\,\cdot\,)
        -(\mathcal{H}V_{0})(0)\,\widetilde{h}_{t}\|_{L^{1}(S)}\,
    <\,\tfrac{\varepsilon}{3}
\end{align}
according to \cref{S1 Main thm 2}.
Under the bi-$K$-invariance assumption, we deduce from \eqref{S4 radial equiv L1}
and \eqref{S4 radial mass} that
\begin{align*}
    (\mathcal{H}V_{0})(0)-(\mathcal{H}v_{0})(0)\,
    =\,\int_{G}\diff{g}\,({V}_{0}(g)-{v}_{0}(g))\,
        \varphi_{0}(g)
    =\,\int_{S}\textrm{d}_{r}{g}\,\widetilde{\delta}(g)^{\frac12}\,
        ({V}_{0}(g)-{v}_{0}(g)).
\end{align*}
Hence, we have thirdly
\begin{align}\label{S4 density3}
    \|((\mathcal{H}V_{0})(0)-(\mathcal{H}v_{0})(0))\widetilde{h}_{t}
    \|_{L^{1}(S)}\,
    \le\,\|\widetilde{v}_{0}-\widetilde{V}_{0}\|_{L^{1}(S)}\,
        \|\widetilde{h}_{t}\|_{L^{1}(S)}\,
    <\,\tfrac{\varepsilon}{3}.
\end{align}
In conclusion, by putting \eqref{S4 density1}, \eqref{S4 density2} and
\eqref{S4 density3} altogether, we obtain
\begin{align*}
    \|\widetilde{v}(t,\,\cdot\,)
        -\widetilde{M}\,\widetilde{h}_{t}\|_{L^{1}(S)}\,
    <\,\varepsilon
\end{align*}
for all $\varepsilon>0$.
Let us turn to the $L^{\infty}$ convergence.
According to \eqref{S4 htilde estimate} and to Theorem 1.6, we have this time
\begin{align*}
    t^{\frac{\ell+|\Sigma_{r}^{+}|}{2}}\,
    \|\widetilde{v}(t,\,\cdot\,)-\widetilde{V}(t,\,\cdot\,)
    \|_{L^{\infty}(S)}\,
    \le\,\|\widetilde{v}_{0}-\widetilde{V}_{0}\|_{L^{1}(S)}\,
        t^{\frac{\ell+|\Sigma_{r}^{+}|}{2}}\,
        \|\widetilde{h}_{t}\|_{L^{\infty}(S)}\,  
    <\,\tfrac{\varepsilon}{3},
\end{align*}
\begin{align*}
    t^{\frac{\ell+|\Sigma_{r}^{+}|}{2}}\,
    \|\widetilde{V}(t,\,\cdot\,)
        -(\mathcal{H}V_{0})(0)\,\widetilde{h}_{t}\|_{L^{\infty}(S)}\,
    <\,\tfrac{\varepsilon}{3},
\end{align*}
and
\begin{align*}
    t^{\frac{\ell+|\Sigma_{r}^{+}|}{2}}\,
    \|((\mathcal{H}V_{0})(0)-(\mathcal{H}v_{0})(0))\widetilde{h}_{t}
    \|_{L^{\infty}(S)}\,
    \le\,
    \|\widetilde{v}_{0}-\widetilde{V}_{0}\|_{L^{1}(S)}\,
        t^{\frac{\ell+|\Sigma_{r}^{+}|}{2}}\,
        \|\widetilde{h}_{t}\|_{L^{\infty}(S)}\,
    <\,\tfrac{\varepsilon}{3}.
\end{align*}
Altogether,
\begin{align*}
    t^{\frac{\ell+|\Sigma_{r}^{+}|}{2}}\,
    \|\widetilde{v}(t,\,\cdot\,)
        -\widetilde{M}\,\widetilde{h}_{t}\|_{L^{\infty}(S)}\,
    \longrightarrow\,0
\end{align*}
as $t\rightarrow\infty$.
The $L^{p}$ convergence follows from convexity.
\end{proof}

\begin{corollary}
The asymptotic convergences \eqref{S1 L1 disting} and \eqref{S1 Linf disting}, 
hence \eqref{S4 Lp disting}, 
still hold with no bi-$K$-invariance condition but under the assumption
\begin{align}\label{S4 other data assumption}
    \int_{G}\diff{g}\,|v_{0}(gK)|e^{\langle{\rho,g^{+}}\rangle}\,
    <\,\infty.
\end{align}
\end{corollary}

\begin{proof}
Notice first that 
\begin{align*}
    \int_{S}\textrm{d}_{r}{g}\,|v_{0}(g)|\widetilde{\delta}(g)^{\frac12}\,
    =\,\int_{G}\diff{g}\,|v_{0}(gK)|\,e^{\langle{\rho,A(g)}\rangle}
    \le\,\int_{G}\diff{g}\,|v_{0}(gK)|e^{\langle{\rho,g^{+}}\rangle}\,
\end{align*}
according to \eqref{S4 ineq Iwasawa Cartan}. 
Hence, the assumption \eqref{S4 other data assumption} is indeed stronger than
$\widetilde{v}_{0}\in{L^{1}(S)}$. Under this assumption, the mass function is
bounded:
\begin{align*}
    |\widetilde{M}(g)|\,
    &\le\,\tfrac{1}{\varphi_{0}(y)}\,
        \int_{G}\diff{g}\,|v_{0}(yK)|\,\varphi_{0}(y^{-1}g)\\[5pt]
    &=\,\tfrac{1}{\varphi_{0}(g)}\,
        \int_{G}\diff{y}\,|v_{0}(yK)|\,\int_{K}\diff{k}\,
        e^{\langle{\rho,A(kg)}\rangle}e^{\langle{\rho,A(ky)}\rangle}\\[5pt]
    &\le\,
        \underbrace{\vphantom{\Big|}
        \tfrac{1}{\varphi_{0}(g)}\,
        \int_{K}\diff{k}\,e^{\langle{\rho,A(kg)}\rangle}
        }_{=\,1}
        \underbrace{\vphantom{\Big|}
        \int_{G}\diff{y}\,|v_{0}(yK)|\,e^{\langle{\rho,y^{+}}\rangle}
        }_{=\,C}\,
    =\,C.
\end{align*}
Here, we used the inequality \eqref{S4 ineq Iwasawa Cartan} and
the fact that $(ky)^{+}=y^{+}$ for all $k\in{K}$.
For proving the $L^{1}$ and the $L^{\infty}$ convergences, 
we argue again by density.
Since $v_{0}(gK)e^{\langle{\rho,g^{+}}\rangle}$ belongs to $L^{1}(\mathbb{X})$,
there exists a function $V_{0}(gK)e^{\langle{\rho,g^{+}}\rangle}$ in
$\mathcal{C}_{c}(\mathbb{X})$ such that 
\begin{align*}
    \int_{G}\diff{g}\,
    |v_{0}(gK)-V_{0}(gK)|\,e^{\langle{\rho,g^{+}}\rangle}\,
    <\,\tfrac{\varepsilon}{3}
\end{align*}
for every $\varepsilon>0$. This implies that the function
$\widetilde{V}_{0}=\widetilde{\delta}^{\tfrac12}V_{0}$ 
approximates the initial data
$\widetilde{v}_{0}=\widetilde{\delta}^{\tfrac12}v_{0}$ in $L^{1}(S)$. Indeed,
\begin{align*}
    \int_{S}\mathrm{d}_{r}{g}\,
    |\widetilde{v}_{0}(g)-\widetilde{V}_{0}(g)|\,
    &=\,\int_{G}\diff{g}\,
        |\widetilde{v}_{0}(gK)-\widetilde{V}_{0}(gK)|\,
        e^{\langle{\rho,A(g)}\rangle}\\[5pt]
    &\le\,\int_{G}\diff{g}\,
        |v_{0}(gK)-V_{0}(gK)|\,
        e^{\langle{\rho,g^{+}}\rangle}\,
    <\,\tfrac{\varepsilon}{3}.
\end{align*}
Let $\widetilde{V}=\widetilde{V}_{0}*\widetilde{h}_{t}$ be the corresponding
solution to the distinguished heat equation, and denote by 
$\widetilde{M}_{V}(g)=\tfrac{(V_{0}*\varphi_{0})(gK)}{\varphi_{0}(g)}$ 
the corresponding mass of $V_{0}$. On the one hand, there exists $T>0$
such that for all $t>T$, we have
\begin{align*}\label{S4 density2}
    \|\widetilde{V}(t,\,\cdot\,)
        -\widetilde{M}_{V}\,\widetilde{h}_{t}\|_{L^{1}(S)}\,
    <\,\tfrac{\varepsilon}{3}
\end{align*}
since $\widetilde{V}_{0}\in\mathcal{C}_{c}(S)$. 
On the other hand, for every $g\in{G}$, we have
\begin{align*}
    |\widetilde{M}(g)-\widetilde{M}_{V}(g)|\,
    &\le\,\tfrac{1}{\varphi_{0}(y)}\,
        \int_{G}\diff{g}\,|v_{0}(yK)-V_{0}(gK)|\,\varphi_{0}(y^{-1}g)\\[5pt]
    &\le\,\int_{G}\diff{y}\,|v_{0}(yK)-V_{0}(gK)|\,
        e^{\langle{\rho,y^{+}}\rangle}\,
    <\,\tfrac{\varepsilon}{3}.
\end{align*}
We conclude by resuming the proof of the previous corollary.
\end{proof}


\vspace{10pt}\noindent\textbf{Acknowledgements.}
The authors are grateful to the referees for checking this manuscript carefully and making helpful suggestions of improvement.
The second author is supported by the Hellenic Foundation for Research and Innovation, Project HFRI-FM17-1733. She acknowledges the SSHN fellowship by the French Institute of Greece and the French Embassy in Greece, which allowed a visit to Institut Denis Poisson, where this work was initiated. Finally, she is thankful to the Institut for the warm hospitality, as well to 
M. Kolountzakis for its help during this stay. 
The last author acknowledges financial support from the Methusalem Programme 
\textit{Analysis and Partial Differential Equations (Grant number 01M01021)} 
during his postdoc stay at Ghent University.
\printbibliography

\vspace{20pt}
\address{
\noindent\textsc{Jean-Philippe Anker:}
\href{mailto:anker@univ-orleans.fr}
{anker@univ-orleans.fr}\\
Institut Denis Poisson,
Universit\'e d'Orl\'eans, Universit\'e de Tours \& CNRS,
Orl\'eans, France}

\vspace{10pt}
\address{
\noindent\textsc{Effie Papageorgiou:}
\href{mailto:papageoeffie@gmail.com}
{papageoeffie@gmail.com}\\
Department of Mathematics and Applied Mathematics,
University of Crete,
Crete, Greece}

\vspace{10pt}
\address{
\noindent\textsc{Hong-Wei Zhang:}
\href{mailto:hongwei.zhang@ugent.be}
{hongwei.zhang@ugent.be}\\
Department of Mathematics:
Analysis, Logic and Discrete Mathematics\\
Ghent University, Ghent, Belgium}
\end{document}

%% file: Concentration1.tex
\begin{tikzpicture}[line cap=round,line join=round,>=triangle 45,x=1cm,y=1cm,scale=2.5]
\clip(-0.5,-1.25) rectangle (5,1.25);
\draw [line width=0.0pt,color=blue,fill=blue,fill opacity=0.31] (0,0) circle (3cm);
\fill[line width=0pt,color=white,fill=white,fill opacity=1] (0,0) -- (4.0654939275664805,0.8360139301518585) -- (0,4) -- (-4,1) -- (-4,-1) -- (0,-4) -- (4.183307817025808,-0.8586935566861592) -- cycle;
\draw [line width=0.5pt] (0,0) circle (2cm);
\draw [line width=0.5pt] (0,0) circle (3cm); 
\fill[line width=0pt,color=white,fill=white,fill opacity=1] (0,0) -- (4,1.5) -- (0,4) -- (-4,1) -- (-4,-1) -- (0,-4) -- (4,-1.5) -- cycle;
 
\draw [line width=0pt,color=white,fill=white,fill opacity=1] (0,0) circle (1.99cm);

\draw [line width=0.5pt] (0,0)-- (4.0654939275664805,0.8360139301518585);
\draw [line width=0.5pt] (0,0)-- (4.0654939275664805,-0.8360139301518585);
\draw [line width=0.5pt] (0,0)-- (4,0);
\draw [line width=0.5pt] (2.5,0) circle (0.5cm);
\draw [color=white,fill=white,fill opacity=1] (1.5,-1) rectangle (2,-1);
\draw [color=white,fill=white,fill opacity=1] (1.5,1) rectangle (2,1);
\draw [color=white,fill=white,fill opacity=1] (2.5,-1) rectangle (3,-1.5);
\draw [color=white,fill=white,fill opacity=1] (2.5,1) rectangle (3,1.5);
\draw [line width=0.5pt] (2.5,0)-- (2.39,0.49);
\draw [line width=0.5pt] (0.75,0) arc (0:30:0.29);
\draw [dashed, line width=0.5pt] (0,0)-- (4.0654939275664805,-1.25);
\draw [dashed, line width=0.5pt] (0,0)-- (4.0654939275664805,1.25);
\draw[decoration={brace,mirror,raise=5pt},decorate]
  (4.2,-1.25) -- node[right=6pt] {$\,\,\overline{\mathfrak{a}^{+}}$} (4.2,1.25);

\begin{scriptsize}
\draw (2.5,0) node{\textbullet};
\draw (2.5,-0.1) node{$2t\rho$};
\draw (2.6,0.25) node{$r(t)$};
\draw (2.2,0.1) node{$\Omega_{t}$};
\draw (1.75,-0.8) node{$|H|=2|\rho|t-r(t)$};
\draw (2.75,-1.1) node{$|H|=2|\rho|t+r(t)$};
\draw (4,-0.1) node{$\rho$-axis};
\draw (0.9,0.09) node{$\gamma(t)$};
\draw (4,1.15) node{wall};
\end{scriptsize}
\end{tikzpicture}

%% file: Concentration2.tex
\begin{tikzpicture}[line cap=round,line join=round,>=triangle 45,x=1cm,y=1cm,scale=2.5]
\clip(-0.5,-1.25) rectangle (5,1.25);
\draw [line width=0.0pt,color=blue,fill=blue,fill opacity=0.31] (0,0) circle (3cm);
\fill[line width=0pt,color=white,fill=white,fill opacity=1] (0,0) -- (4.0654939275664805,0.8360139301518585) -- (0,4) -- (-4,1) -- (-4,-1) -- (0,-4) -- (4.183307817025808,-0.8586935566861592) -- cycle;
\draw [line width=1.5pt] (0,0) circle (1cm);
\draw [line width=0.5pt] (0,0) circle (3cm);   
\fill[line width=0pt,color=white,fill=white,fill opacity=1] (0,0) -- (4,1.5) -- (0,4) -- (-4,2) -- (-4,-2) -- (0,-4) -- (4,-1.5) -- cycle;
\draw [line width=0pt,color=white,fill=white,fill opacity=1] (0,0) circle (1cm);

\draw [line width=0.5pt] (0,0)-- (4.0654939275664805,0.8360139301518585);
\draw [line width=0.5pt] (0,0)-- (4.0654939275664805,-0.8360139301518585);
\draw [color=white,fill=white,fill opacity=1] (1.5,-1) rectangle (2,-1);
\draw [color=white,fill=white,fill opacity=1] (1.5,1) rectangle (2,1);
\draw [color=white,fill=white,fill opacity=1] (2.5,-1) rectangle (3,-1.5);
\draw [color=white,fill=white,fill opacity=1] (2.5,1) rectangle (3,1.5);
\draw [dashed, line width=0.5pt] (0,0)-- (4.0654939275664805,-1.25);
\draw [dashed, line width=0.5pt] (0,0)-- (4.0654939275664805,1.25);
\draw[decoration={brace,mirror,raise=5pt},decorate]
  (4.2,-1.25) -- node[right=6pt] {$\,\,\overline{\mathfrak{a}^{+}}$} (4.2,1.25);

\begin{scriptsize}
\draw (0.9,-0.5) node{$|H|=\varepsilon(t)\sqrt{t}$};
\draw (2.8,-1.1) node{$|H|=\tfrac{\sqrt{t}}{\varepsilon(t)}$};
\draw (2,0) node{$\widetilde{\Omega}_{t}$};
\draw (4,1.15) node{wall};
\end{scriptsize}
\end{tikzpicture}